\newcommand\CC{{\mathbb C}}
\newcommand\cA{{\cal A}}

\newcommand\cC{{\cal C}}
\newcommand\cD{{\cal D}}

\newcommand\cI{{\cal I}}

\newcommand\cN{{\cal N}}
\newcommand{\contr}{{\mspace{1.5mu}\lrcorner\mspace{2mu}}}
\newcommand\cO{{\cal O}}

\newcommand\cS{{\cal S}}

\newcommand\cU{{\cal U}}
\newcommand\cV{{\cal V}} 
\newcommand\cW{{\cal W}}
\newcommand\cX{{\cal X}} 
\newcommand\cY{{\cal Y}}

\newcommand\dra{\dashrightarrow}
\newcommand\es{\emptyset}
\newcommand\Ext{\text{Ext}}

\newcommand\hra{\hookrightarrow}
\newcommand\la{\langle}
\newcommand\lag{\mathbb{LG}}

\newcommand\lagr{\mathbb{LG}(\bigwedge^ 3 V)}

\newcommand\lra{\longrightarrow}
\newcommand\n{\noindent}
\newcommand\NN{{\mathbb N}}
\newcommand\ov{\overline}
\newcommand\PP{{\mathbb P}}

\newcommand\ra{\rangle}

\newcommand\wh{\widehat}
\newcommand\wt{\widetilde}
\newcommand\ZZ{{\mathbb Z}}
\newcommand{\Gr}{\mathrm{Gr}}
\newcommand{\mapor}[1]{{\stackrel{#1}{\longrightarrow}}}
\newcommand{\mapver}[1]{\Big\downarrow
\vcenter{\rlap{$\scriptstyle#1$}}}
\documentclass[10pt,a4paper]{article}
\usepackage{amsthm}
\usepackage{amsmath}
\usepackage{amssymb}
\usepackage{makeidx}         
\usepackage{graphicx}        
\usepackage{multicol}        
\usepackage[all,ps]{xy}
\usepackage{layout}
\theoremstyle{plain}

\newtheorem{thm}{Theorem}[section]
\newtheorem{clm}[thm]{Claim}

\newtheorem{crl}[thm]{Corollary}
\newtheorem{hyp}[thm]{Hypothesis}
\newtheorem{lmm}[thm]{Lemma}
\newtheorem{prp}[thm]{Proposition}
\newtheorem{prp-dfn}[thm]{Proposition-Definition}

\theoremstyle{definition}

\newtheorem{ass}[thm]{Assumption}
\newtheorem{dfn}[thm]{Definition}

\theoremstyle{remark}

\newtheorem{expl}[thm]{Example}

\newtheorem{rmk}[thm]{Remark}

\DeclareMathOperator{\Ann}{Ann}

\DeclareMathOperator{\cod}{cod}

\DeclareMathOperator{\coker}{coker}
\DeclareMathOperator{\cork}{cork}
\DeclareMathOperator{\Def}{Def}

\DeclareMathOperator{\Fix}{Fix}

\DeclareMathOperator{\GL}{Gl}
\DeclareMathOperator{\Hilb}{Hilb}
\DeclareMathOperator{\Hom}{Hom}
\DeclareMathOperator{\Id}{Id}
\DeclareMathOperator{\im}{im}

\DeclareMathOperator{\Pic}{Pic}

\DeclareMathOperator{\PGL}{PGL}

\DeclareMathOperator{\rk}{rk}

\DeclareMathOperator{\sing}{sing}
\DeclareMathOperator{\Spec}{Spec}

\DeclareMathOperator{\supp}{supp}
\DeclareMathOperator{\Sym}{Sym}
\DeclareMathOperator{\vol}{vol}

\usepackage{ifthen}
\newcommand{\cit}[1]{{\rm \textbf{#1}}}
\newcommand{\Ref}[2]{\cit{%
\ifthenelse{\equal{#1}{thm}}{Theorem}{}%
\ifthenelse{\equal{#1}{ass}}{Assumption}{}%
\ifthenelse{\equal{#1}{chp}}{Chapter}{}%
\ifthenelse{\equal{#1}{prp}}{Proposition}{}%
\ifthenelse{\equal{#1}{lmm}}{Lemma}{}%
\ifthenelse{\equal{#1}{crl}}{Corollary}{}%
\ifthenelse{\equal{#1}{dfn}}{Definition}{}%
\ifthenelse{\equal{#1}{expl}}{Example}{}%
\ifthenelse{\equal{#1}{hyp}}{Hypothesis}{}%
\ifthenelse{\equal{#1}{rmk}}{Remark}{}%
\ifthenelse{\equal{#1}{clm}}{Claim}{}%
\ifthenelse{\equal{#1}{exe}}{Exercise}{}%
\ifthenelse{\equal{#1}{sec}}{Section}{}%
\ifthenelse{\equal{#1}{subsec}}{Subsection}{}%
\ifthenelse{\equal{#1}{univ}}{Universal Property}{}%
\ifthenelse{\equal{#1}{trm}}{Terminology}{}%
\ifthenelse{\equal{#1}{tbl}}{Table}{}%
\  \ref{#1:#2}%
}}

\usepackage{multirow}

 \makeindex
 
\begin{document}
 \title{Double covers of EPW-sextics}
 \author{Kieran G. O'Grady\thanks{Supported by
 PRIN 2007}\\\\
\lq\lq Sapienza\rq\rq Universit\`a di Roma}
\date{June 8 2012}
 \maketitle
 \tableofcontents
 \section{Introduction}\label{sec:prologo}
 \setcounter{equation}{0}
EPW-sextics are defined as follows. Let $V$ be a $6$-dimensional complex vector space. Choose a volume-form  
$\vol\colon\bigwedge^6 V\overset{\sim}{\lra}\CC$
 and  equip $\bigwedge^3 V$ with the symplectic form
\begin{equation}
  (\alpha,\beta)_V:=\vol(\alpha\wedge\beta).
\end{equation}
Let $\lagr$ be the symplectic Grassmannian parametrizing Lagrangian subspaces of $\bigwedge^3 V$ - of course $\lagr$ does not depend on the choice of volume-form. 
Let $F\subset\bigwedge^3 V\otimes\cO_{\PP(V)}$ be the sub vector-bundle 
 with fiber  
\begin{equation}
F_v:=\{\alpha\in\bigwedge^3 V\mid v\wedge\alpha=0\}
\end{equation}
 over $[v]\in\PP(V)$.  Notice that $(,)_V$ is zero on $F_v$ and   $2\dim(F_v)=20=\dim\bigwedge^3 V$; thus $F$ is a Lagrangian sub vector-bundle of the trivial symplectic vector-bundle on $\PP(V)$ with fiber $\bigwedge^3 V$. 
Next choose $A\in\lagr$. Let
\begin{equation}\label{diecidieci}
F\overset{\lambda_A}{\lra}(\bigwedge^3 V/A)\otimes\cO_{\PP(V)}
\end{equation}
be the composition of the inclusion $F\subset\bigwedge^3 V\otimes\cO_{\PP(V)}$ followed by the quotient map.
Since $\rk F=\dim(V/A)$ 
 the determinat of $\lambda_A$ makes sense. 
    Let 
   \begin{equation*}
Y_A:=V(\det\lambda_A).
\end{equation*}
A straightforward computation gives that $\det F\cong\cO_{\PP(V)}(-6)$ and hence 
   $\det\lambda_A\in H^0(\cO_{\PP(V)}(6))$. 
It follows that if   $\det\lambda_A\not=0$ then  $Y_A$ is  a sextic hypersurface. As is easily checked $\det\lambda_A\not=0$ for  generic $A\in\lagr$  (notice that there exist \lq\lq pathological\rq\rq $A$'s such that  $\lambda_A=0$ e.g.~$A=F_{v_0}$).
An {\it EPW-sextic} (after Eisenbud, Popescu and Walter~\cite{epw}) is a sextic hypersurface in $\PP^5$ which is projectively equivalent to $Y_A$ for some $A\in\lagr$.
Let $Y_A$ be an  EPW-sextic. 
One constructs a coherent sheaf $\xi_A$ on $Y_A$ and a multiplication map $\xi_A\times\xi_A\to \cO_{Y_A}$ which gives $\cO_{Y_A}\oplus\xi_A$ a structure of $\cO_{Y_A}$-algebra - this is known to experts, see~\cite{catanese} - we will give the construction in~\Ref{subsec}{voiladouble}. The \emph{double EPW-sextic} associated to $A$ is $X_A:=\Spec(\cO_{Y_A}\oplus\xi_A)$; we let $f_A\colon X_A\to Y_A$ be the structure morphism. In~\cite{og2} we considered $X_A$ for generic $A$ and we proved that it is a Hyperk\"ahler deformation of  $(K3)^{[2]}$ (the blow-up of the diagonal in the symmetric square of a $K3$ surface). In the present paper we will analyze $X_A$ for $A$ varying in a codimension-$1$ subset of $\lagr$. In order to state our main results we will introduce some notation. Given $A\in\lagr$  we let 
\begin{eqnarray}
Y_A(k)= & \{[v]\in\PP(V)\mid \dim(A\cap F_v)= k\}, \\
Y_A[k]= & \{[v]\in\PP(V)\mid \dim(A\cap F_v)\ge k\}. 
\end{eqnarray}
Thus $Y_A(0)=(\PP(V)\setminus Y_A)$ and $Y_A=Y_A[1]$. 
Double EPW-sextics come with a natural polarization; we let 
\begin{equation}\label{ecceampio}
\cO_{X_A}(n):=f_A^{*}\cO_{Y_A}(n),\quad
H_A\in |\cO_{X_A}(1)|.
\end{equation}
 The following closed subsets of $\lagr$ play a key r\^ole in the present paper:
\begin{eqnarray}
\Sigma:= & \{A\in\lagr\mid \text{$\exists W\in{\mathbb G}r(3,V)$ s.~t.~$\bigwedge^3 W\subset A$}\},\\
\Delta:=  &  \{A\in\lagr\mid Y_A[3]\not=\es\}.
\end{eqnarray}
A straightforward computation,  see~\cite{og5}, gives that $\Sigma$ is  irreducible of codimension $1$. A similar computation, see~\Ref{prp}{codelta}, gives that $\Delta$  is  irreducible of codimension $1$ and distinct from $\Sigma$.
Let
\begin{equation}\label{eccozero}
\lagr^0:=  \lagr\setminus\Sigma\setminus\Delta\,.
\end{equation}
Thus $\lagr^0$ is open dense in $\lagr$. 
In~\cite{og2} we proved that if $A\in\lagr^0$ then  $X_A$ is a hyperk\"ahler (HK) $4$-fold which can be deformed to $(K3)^{[2]}$.
Moreover we showed that the family of polarized HK $4$-folds $ (X_A,H_A)$ for $A$ varying in$\lagr^0$
 is  locally complete.  Three other explicit locally complete families of projective HK's of dimension greater than $2$ are known - see~\cite{beaudon,debvoi,iliran1,iliran2}. In all of the  examples the HK manifolds are deformations of  the Hilbert square of a $K3$: they are distinguished by the value of the Beauville-Bogomolov form on the polarization class (it equals $2$ in the case of double EPW-sextics and $6$, $22$ and $38$ in the other cases). 
  In the present paper we will analyze $X_A$ for $A\in\Delta$, mainly under the hypothesis that $A\not\in\Sigma$. 
Let $A\in(\Delta\setminus\Sigma)$. We will prove the following results
 \begin{enumerate}
\item[(1)]
$Y_A[3]$ is a finite set and it equals $Y_A(3)$.  If $A$ is generic in $(\Delta\setminus\Sigma)$  then $Y_A(3)$ is  a singleton.
\item[(2)]
One may associate to $[v_0]\in Y_A(3)$ a $K3$ surface $S_A(v_0)\subset\PP^6$ of genus $6$, well-defined up to projectivities. Conversely the generic $K3$ of genus $6$ is projectively equivalent to $S_A(v_0)$ for some $A\in(\Delta\setminus\Sigma)$ and $[v_0]\in Y_A(3)$.
\item[(3)]
The singular set of $X_A$ is equal to $f_A^{-1}Y_A(3)$. There is a single $p_i\in X_A$ mapping to  $[v_i]\in Y_A(3)$ and the cone of $X_A$ at $p_i$ is isomorphic to the cone over the set of incident couples $(x,r)\in\PP^2\times(\PP^2)^{\vee}$  (i.e.~$\PP(\Omega_{\PP^2})$). Thus we have two  standard small resolutions of a neighborhood of $p_i$ in $X_A$, one with fiber $\PP^2$ over $p_i$, the other with fiber $(\PP^2)^{\vee}$. Making a choice $\epsilon$ of local small resolution at each $p_i$ we get 
 a resolution $X^{\epsilon}_A\to X_A$ 
with the following properties: There is a birational map $X^{\epsilon}_A\dashrightarrow S_A(v_i)^{[2]}$ such that the pull-back of a holomorphic symplectic form on $S_A(v_i)^{[2]}$ is a symplectic form on $X^{\epsilon}_A$.  If  $S_A(v_i)$ contains no lines (true for generic $A$ by Item~(2)) then there exists a choice of $\epsilon$ such that $X^{\epsilon}_A$  is isomorphic to $S_A(v_i)^{[2]}$. 
\item[(4)]
Given a sufficiently small open (classical topology) $\cU\subset(\lagr\setminus\Sigma)$ containing $A$ the family of double EPW-sextics parametrized by $\cU$ has a simultaneous resolution of singularities (no base change) with fiber $X_A^{\epsilon}$ over $A$ (for an arbitrary choice of $\epsilon$).
\end{enumerate}
A remark: if $Y_A(3)$ has more than one point we do not expect all the small resolutions to be projective (i.e.~K\"ahler).
 Items~(1)-(4) should be compared with known results on cubic $4$-folds - recall that if $Z\subset\PP^5$ is a smooth cubic hypersurface  the variety $F(Z)$ parametrizing lines in $Z$ is a HK $4$-fold which can be deformed to $(K3)^{[2]}$ and moreover the primitive weight-$4$ Hodge structure of $Z$ is isomorphic (after a Tate twist) to the  primitive weight-$2$ Hodge structure of $F(Z)$, see~\cite{beaudon}. Let $D\subset|\cO_{\PP^5}(3)|$ be the prime divisor parametrizing singular cubics. Let $Z\in D$ be generic: the following results are well-known.
 \begin{enumerate}
\item[(1')]
 $\sing Z$ is a finite set.
\item[(2')]
Given $p\in\sing Z$  the set  $S_Z(p)\subset F(Z)$ of lines containing $p$ is a $K3$ surface of genus $4$ and viceversa the generic such $K3$ is isomorphic to $S_Z(p)$ for some $Z$ and $p\in\sing Z$. 
\item[(3')]
$F(Z)$ is birational to $S_Z(p)^{[2]}$.
\item[(4')]
After a local base-change of order $2$ ramified along $D$ the period map extends across $Z$.
\end{enumerate}
Thus Items~(1')-(2')-(3') are analogous to Items~(1), (2) and~(3) above, Item~(4') is analogous to (4) but there is  an important difference namely the need for a base-change of order $2$. (Actually the paper~\cite{og3} contains results showing that there is a statement valid for cubic hypersurfaces which is even closer to our result for double EPW-sextics, 
the r\^ole of $\Sigma$ being played by the divisor parametrizing cubics containing a plane.)  
 We explain the relevance of Items~(1)-(4). Items~(3) and~(4)  prove  the theorem of ours mentioned above i.e.~that if $A\in\lagr^0$ then $X_A$ is a HK deformation of $(K3)^{[2]}$ (the family of polarized double EPW-sextics is locally complete by a straightforward parameter count). The proof in this paper is independent of the proof in~\cite{og2}. Beyond giving a new proof of an \lq\lq old\rq\rq theorem the above results show that away from $\Sigma$ the period map is regular,  it lifts (locally) to the relevant classifying space and the value at $A\in(\Delta\setminus\Sigma)$ may be identified with the period point of the Hilbert square $S_A(v_0)^{[2]}$. We remark that in~\cite{og4} we had proved that the period map is as well-behaved as possible at the generic $A\in(\Delta\setminus\Sigma)$, however we did not have the exact statement about $X_A^{\epsilon}$ and we had no statement about an arbitrary $A\in(\Delta\setminus\Sigma)$. 

The paper is organized as follows. In~\Ref{sec}{epwdoppie}  we will give formulae that  describe double EPW-sextics locally. The formulae are  known to experts, see~\cite{catanese}, we will go through the proofs because we could not find a suitable reference. We will also perform the local computations needed to prove  Item~(4) above. 
 In~\Ref{sec}{avoidelta} we will go through some standard computations involving $\Delta$. 
In~\Ref{sec}{famdesing} we will prove Items~(1), (4) and the statements of Item~(3) which do not involve the $K3$ surface $S_A(v_0)$. In~\Ref{sec}{zitelle} we will prove Item~(2) and the remaining statement of Item~(3).  \Ref{sec}{wedding} contains auxiliary results on $3$-dimensional linear sections of $\Gr(3,\CC^5)$. 
\vskip 3mm
\n
{\bf Notation and conventions:} 
Throughout the paper $V$ is a $6$-dimensional complex vector space. 
\vskip 2mm
\noindent
Let $W$ be a finite-dimensional complex vector-space. The span of a subset $S\subset W$ is denoted by $\la S\ra$. Let $S\subset\bigwedge^q W$.  The {\it support of $S$} is the smallest subspace $U\subset W$ such that $S\subset\im(\bigwedge^{q}U\lra \bigwedge^q W)$: we denote it by $\supp(S)$, if $S=\{\alpha\}$ is a singleton we let $\supp(\alpha)=\supp(\{\alpha\})$  (thus if $q=1$ we have $\supp(\alpha)=\la\alpha\ra$).  We define the support of a set of symmetric tensors analogously. If  $\alpha\in\bigwedge^q W$ or $\alpha\in \Sym^d W$ the {\it rank of $\alpha$} is the  dimension of $\supp(\alpha)$. 
An element of $\Sym^2 W^{\vee}$ may be viewed either  
as a symmetric map or as a quadratic form: we will denote the former by $\wt{q},\wt{r},\ldots$ and the latter by $q,r,\ldots$ respectively.
\vskip 2mm
\noindent
Let $M=(M_{ij})$ be a $d\times d$ matrix with entries in a commutative ring $R$. We let  $M^c=(M^{ij})$ be the matrix of cofactors of $M$, i.e.~$M^{i,j}$ is $(-1)^{i+j}$ times the determinant of the matrix obtained from $M$    
 by deleting its $j$-th row and $i$-th column. 
 We recall the following  interpretation of $M^c$. Suppose that $f\colon A\to B$ is a linear map between free $R$-modules of rank $d$ and that $M$ is the matrix associated to $f$ by the  choice of bases $\{a_1,\ldots,a_d\}$ and $\{b_1,\ldots,b_d\}$ of $A$ and $B$ respectively. Then $\bigwedge^{d-1}f$ may be viewed as a map
 \begin{equation}
\bigwedge^ {d-1}f\colon A^{\vee}\otimes\bigwedge^d A\cong
\bigwedge^{d-1} A \lra
\bigwedge^{d-1} B\cong B^{\vee}\otimes\bigwedge^d B.
\end{equation}
(Here $A^{\vee}:=\Hom(A,R)$ and similarly for $B^{\vee}$.)
The matrix associated to $\bigwedge^ {d-1}f$ by the choice of bases $\{a^{\vee}_1\otimes (a_1\wedge\ldots\wedge a_d),\ldots,a^{\vee}_d\otimes (a_1\wedge\ldots\wedge a_d)\}$ and $\{b^{\vee}_1\otimes (b_1\wedge\ldots\wedge b_d),\ldots,b^{\vee}_d\otimes (b_1\wedge\ldots\wedge b_d)\}$ is equal to $M^c$. 
\vskip 2mm
\noindent
Let $W$ be a finite-dimensional complex vector-space. We will adhere to pre-Grothendieck conventions: $\PP(W)$ is the set of $1$-dimensional vector subspaces of $W$. Given  a non-zero $w\in W$ 
we will denote the span of $w$ by $[w]$ rather than $\la w\ra$; this agrees with  standard notation. 
Suppose that $T\subset\PP(W)$. Then $\la T\ra\subset\PP(W)$ is the {\it projective span of  $T$} i.e.~the intersection of all linear subspaces of $\PP(W)$ containing $T$. 
\vskip 2mm
\noindent
Schemes  are defined over $\CC$, the topology is the Zariski topology unless we state the contrary. Let $W$ be finite-dimensional complex vector-space: $\cO_{\PP(W)}(1)$ is the line-bundle on $\PP(W)$ with fiber $L^{\vee}$ on the point $L\in\PP(W)$. Let $F\in \Sym^d W^{\vee}$: we let  $V(F)\subset\PP(W)$ be the subscheme defined by vanishing of $F$. If $E\to X$ is a vector-bundle we denote by $\PP(E)$ the projective fiber-bundle with fiber $\PP(E(x))$ over $x$ and we define $\cO_{\PP(W)}(1)$ accordingly.
If $Y$ is a subscheme of $X$ we let $Bl_Y X\lra X$ be the blow-up of $Y$.
 \section{Symmetric resolutions and double covers}\label{sec:epwdoppie}
 \setcounter{equation}{0}  
In~\Ref{subsec}{prodenne} we will describe a method (well-known to experts) for constructing double covers. In~\Ref{subsec}{voiladouble} we will show how to implement the construction in order to construct double EPW-sextics.  \Ref{subsec}{doppipiccoli} contains the main ingredients needed to construct the simultaneous desingularization described in Item~(3) of~\Ref{sec}{prologo}. 
\subsection{Product formula and double covers}\label{subsec:prodenne}
\setcounter{equation}{0}
Let $R$ be an integral  Noetherian ring. 
Let $N$ be an  $R$-module with a free resolution
\begin{equation}\label{risenne}
0\lra U_1\overset{\lambda}{\lra} U_0\overset{\pi}{\lra} N\lra 0,\qquad \rk\, U_1=\rk\, U_0=d>0.
\end{equation}
Let $\{a_1,\ldots,a_d\}$ and $\{b_1,\ldots,b_d\}$  be  bases of $U_0$ and $U_1$ respectively. 
  Let $M_{\lambda}$  be the matrix associated to $\lambda$  by our choice of bases - notice that $\det M_{\lambda}$ annihilates $N$. 
  Given  a homomorphism 
\begin{equation}\label{morbeta}
\beta\colon N\to \Ext^1(N,R) 
\end{equation}
one defines a product $m_{\beta}\colon N\times N\to R/(\det M_{\lambda})$
as follows. 
Applying the $\Hom(\,\cdot\,,R)$-functor to~\eqref{risenne} we get the exact sequence
\begin{equation}\label{risext}
0\lra U^{\vee}_0\overset{\lambda^t}{\lra} U^{\vee}_1\overset{\rho}\lra \Ext^1(N,R)\lra 0.
\end{equation}
In particular $\det M_{\lambda}$ kills $\Ext^1(N,R)$. Now apply
 the functor $\Hom(N,\,\cdot\,)$ to the  exact sequence 
\begin{equation}
0  \lra  R  \overset{\det M_{\lambda}}\lra  R   \lra   
 R/ (\det M_{\lambda}) \lra  0. \\
\end{equation}
Since $\Ext^1(N,R)\to \Ext^1(N,R)$ is multiplication by $\det M_{\lambda}$ we get 
 a coboundary isomorphism
\begin{equation}\label{isobordo}
\partial\colon \Hom(N,R/(\det M_{\lambda}))\overset{\sim}{\lra} 
\Ext^1(N,R).
\end{equation}
We let
 \begin{equation}\label{emmebeta}
\begin{matrix}
N\times N & \overset{m_{\beta}}{\lra} &  R/(\det M_{\lambda})\\
 (n,n') & \mapsto & (\partial^{-1}\beta(n))(n').
\end{matrix}
\end{equation}
We will give an explicit formula for $m_{\beta}$. 
Let $\pi\colon U_0\to N$ be as in~\eqref{risenne}. 
Then $\beta\circ\pi$ lifts to a homomorphism $\mu^t\colon U_0\to U_1^{\vee}$ (the map is written as  a transpose in order to conform to the notation for double EPW-sextics - see~\Ref{subsec}{voiladouble}). It follows that there exists $\alpha\colon U_1\to U^{\vee}_0$ such that  
\begin{equation}\label{unozeroenne}
\begin{array}{ccccccccc}
0 & \to &  U_1 &\mapor{\lambda}& U_0 & \overset{\pi}{\lra} & N & \to & 0\\
 & & \mapver{\alpha}& &\mapver{\mu^t} & &
\mapver{\beta}& & \\
0 & \to & U_0^{\vee} & \mapor{\lambda^{t}} & U_1^{\vee} & \overset{\rho}{\lra} & \Ext^1(N,R) & \to & 0
\end{array}
\end{equation}
is  a commutative diagram.
Let $\{a^{\vee}_1,\ldots,a^{\vee}_d\}$ and $\{b^{\vee}_1,\ldots,b^{\vee}_d\}$  be the  bases of $U_0^{\vee}$ and $U_1^{\vee}$ which are dual to the chosen bases of $U_0$ and $U_1$. 
 Let $M_{\mu^t}$ be the matrix associated to $\mu^t$ by our choice of bases. 
\begin{prp}
Keeping notation as above we have
\begin{equation}\label{prodmin}
m_{\beta}(\pi(a_i),\pi(a_j))\equiv (M^c_{\lambda}\cdot M_{\mu^t})_{ji}\mod{(\det M_{\lambda})} 
\end{equation}
where $M^c_{\lambda}$ is the matrix of cofactors of $M_{\lambda}$. 
\end{prp}
 \begin{proof}
Equation~\eqref{risext} gives an isomorphism 
\begin{equation}\label{ecconu}
\nu\colon \Ext^1(N,R)\overset{\sim}{\to} U_1^{\vee}/\lambda^t(U_0^{\vee}).
\end{equation}
Let $\det(U_{\bullet}):=\bigwedge^d U_1^{\vee}\otimes\bigwedge^d U_0$. 
We will define an isomorphism
\begin{equation}\label{tetaisom}
\theta\colon U_1^{\vee}/\lambda^t(U_0^{\vee})\overset{\sim}{\lra} 
\Hom\left(N,\det(U_{\bullet})/(\det\lambda)\right).
\end{equation}
First let
\begin{equation}\label{mappazza}
\begin{matrix}
U_1^{\vee}=& \bigwedge^{d-1}U_1\otimes\bigwedge^d U_1^{\vee} &
\overset{\wh{\theta}}{\lra} & 
\bigwedge^{d-1}U_0\otimes\bigwedge^d U_1^{\vee} & = &
\Hom(U_0,\det (U_{\bullet})) \\
&\zeta\otimes\xi & \mapsto & \bigwedge^ {d-1}(\lambda)(\zeta)\otimes\xi\,.
\end{matrix}
\end{equation}
We claim that 
\begin{equation}\label{jamesbiondi}
\im(\wh{\theta})=\{\phi\in \Hom(U_0,\det(U_{\bullet}))\mid 
\phi\circ\lambda(U_1)\subset(\det\lambda)\}.
\end{equation}
In fact by Cramer's formula 
\begin{equation}\label{cofattori}
M_{\lambda}^c\cdot M_{\lambda}^t=M_{\lambda}^t\cdot M_{\lambda}^c=\det M_{\lambda}\cdot 1
\end{equation}
and Equation~\eqref{jamesbiondi} follows.
 Thus $\wh{\theta}$ induces  a surjective homomorphism
 \begin{equation}
\wt{\theta}\colon U_1^{\vee}\lra 
\Hom\left(N, \det (U_{\bullet})/(\det\lambda)\right).
\end{equation}
One checks easily  that 
$\lambda^t(U_0^{\vee})= \ker\wt{\theta}$ - use Cramer again. 
We define   $\theta$  to be the homomorphism  induced by   $\wt{\theta}$; we have proved that it is an isomorphism.
We claim that
\begin{equation}\label{voila}
\theta\circ\nu=\partial^{-1},\qquad\text{$\partial$ as in~(\ref{isobordo}).}
\end{equation}
In fact let $K$ be the fraction field of $R$ and $0\to R\overset{\iota}{\to} I^0\to I^1\to\ldots$ be an injective resolution of $R$ with $I^0=\det(U_{\bullet})\otimes K$ and $\iota(1)=\det\lambda\otimes 1$.
Then $\Ext^{\bullet}(N,R)$ is the cohomology of the double complex $\Hom(U_{\bullet},I^{\bullet})$ and of course also of the single complexes $\Hom(U_{\bullet},R)$ and $\Hom(N,I^{\bullet})$. 
One checks easily that the isomorphism $\partial$ of~(\ref{isobordo}) is equal to the isomorphism $H^1(\Hom(N,I^{\bullet}))\overset{\sim}{\to}H^1(\Hom(U_{\bullet},I^{\bullet}))$ i.e.
\begin{equation}\label{eccoeps}
\partial\colon
\Hom(N,\det (U_{\bullet})/(\det\lambda))=
\Hom(N,I^0/\iota(R))\overset{\sim}{\lra}H^1(\Hom(U_{\bullet},I^{\bullet})).
\end{equation}
Let
 $f\in \Hom(N,\det (U_{\bullet})/(\det\lambda))$; a representative of $\partial(f)$ in the double complex $\Hom(U_{\bullet},I^{\bullet})$ is given by $g^{0,1}:=f\circ\pi\in \Hom(U_0,I^1)$. Let $g^{0,0}\in \Hom(U_0,\det( U_{\bullet}))$ be a lift of $g^{0,1}$ and $g^{1,0}\in \Hom(U_1,\det( U_{\bullet}))$ be defined by $g^{1,0}:=g^{0,0}\circ\lambda$. One checks that $\im(g^{1,0})\subset(\det\lambda)$ and hence there exists $g\in \Hom(U_1,R)$ such that  $g^{1,0}=\iota\circ g$. By construction $g$ represents a class $[g]\in H^1(\Hom(U_{\bullet},R))=U_1^{\vee}/\lambda^t(U_0^{\vee})$ and $[g]=\nu\circ\partial(f)$. An explicit computation shows that $[g]=\theta^{-1}(f)$. 
 This proves~(\ref{voila}). Now we prove Equation~(\ref{prodmin}).
By~(\ref{voila}) we have
\begin{equation}
m_{\beta}(\pi(a_i),\pi(a_j))=(\partial^{-1}\beta\pi(a_i))(\pi(a_j))
=(\theta\nu\beta\pi(a_i))(\pi(a_j)).
\end{equation}
Unwinding the definition of $\theta$ one gets that the right-hand side of the above equation equals  the right-hand side of~(\ref{prodmin}).
 \end{proof}
Let $m_{\beta}$ be given by~\eqref{emmebeta}: we define   a product on $R/(\det M_{\lambda})\oplus N$ as follows. Let $(r,n),(r',n')\in R/(\det M_{\lambda})\oplus N$: we set
\begin{equation}\label{moltiplico}
(r,n)\cdot(r',n'):=(rr'+m_{\beta}(n,n'),rn'+r'n).
\end{equation}
In general the above product is neither  associative nor commutative. We will give an example in which the product is both associative and commutative. 
Suppose that we have
\begin{equation}\label{mappasimm}
0\lra U^{\vee} \overset{\gamma}{\lra} U \overset{\pi}{\lra} N \lra 0,\qquad \gamma^t=\gamma
\end{equation}
with $U$ a free $R$-module of rank $d>0$ and the  sequence is supposed to be exact.
 We get a commutative diagram~\eqref{unozeroenne} by letting
 \begin{equation*}
U_0:=U,\quad U_1:=U^{\vee},\quad\lambda=\gamma,\quad\alpha=\Id _{U^{\vee}},\quad
\mu^t=\Id _U,
\end{equation*}
 and  $\beta=\beta(\gamma)\colon N\to \Ext^1(N,R)$  the  map  induced by $\Id _U$.
Abusing notation we let $m_{\gamma}\colon N\times N\to R/(\det M_{\gamma})$ be the map defined by~$m_{\beta(\gamma)}$. 
\begin{prp}\label{prp:associo}
Suppose that we have Exact Sequence~\eqref{mappasimm}. The product  on $R/(\det M_{\gamma})\oplus N$ defined  by $m_{\gamma}$  is  associative and commutative.
\end{prp}
\begin{proof}
Let $d:=\rk\, U>0$.
Let $\{a_1,\ldots,a_d\}$ be a basis of $U$ and $\{a^{\vee}_1,\ldots,a^{\vee}_d\}$ be the dual basis of $U^{\vee}$. Let $M=M_{\gamma}$ i.e.~the matrix associated to $\gamma$ by our choice of bases. 
 By~\eqref{prodmin} we have
\begin{equation}\label{mucca}
m_{\gamma}(\pi(a_i),\pi(a_j))\equiv M^c_{ji} \mod{(\det M)}.
\end{equation}
Since $\gamma$ is a symmetric map $M$ is a symmetric matrix. Thus  $M^c$  is a symmetric matrix. By~\eqref{mucca} we get that $m_{\gamma}$ is symmetric. It remains to prove that $m_{\gamma}$ is associative. 
For $1\le i< k\le d$ and $1\le h\not= j\le d$ let $M^{i,k}_{h,j}$ be the $(d-2)\times(d-2)$-matrix obtained by deleting from $M$   rows $i,k$ and columns $h,j$. 
Let  $X_{ijk}=(X^h_{ijk})\in R^d$ be defined by
\begin{equation}
X^h_{ijk}:=
\begin{cases}
(-1)^{i+k+j+h}\det M^{i,k}_{j,h} & \text{if $h< j$}, \\
0 & \text{if $h=j$}. \\
(-1)^{i+k+j+h-1}\det M^{i,k}_{j,h} & \text{if $j<h$.}
\end{cases}
\end{equation}
A tedious but straightforward computation gives that
\begin{equation}\label{relazione}
M_{ij}^c a_k-M^c_{jk} a_i=\gamma(\sum_{h=1}^d X^h_{ijk}a_h^{\vee}).
\end{equation}
The above equation proves associativity of $m_{\gamma}$.  
\end{proof}
Keep hypotheses as in~\Ref{prp}{associo}.  We let
\begin{equation}
X_{\gamma}:= \Spec(R/(\det M_{\lambda})\oplus N),\qquad  
Y_{\gamma}:=\Spec(R/(\det M_{\lambda})).
\end{equation}
Let  $f_{\gamma}\colon X_{\gamma}\to Y_{\gamma}$ be the structure map. 
We realize $X_{\gamma}$ as a subscheme of $\Spec(R[\xi_1,\ldots,\xi_d])$ as follows. Since the ring $R/(\det M_{\gamma})\oplus N$
is associative and commutative there is a well-defined surjective morphism of $R$-algebras
\begin{equation}
R[\xi_1,\ldots,\xi_d]\longrightarrow R/(\det M_{\gamma})\oplus N
\end{equation}
mapping $\xi_i$ to $a_i$. Thus we have an inclusion 
\begin{equation}\label{esplicito}
X_{\gamma}\hra \Spec(R[\xi_1,\ldots,\xi_d]). 
\end{equation}
\begin{clm}\label{clm:ideawu}
Referring to Inclusion~(\ref{esplicito})  the ideal  of $X_{\gamma}$ is generated by the entries of the matrices
\begin{equation}\label{profumo}
M_{\gamma}\cdot\xi,\qquad \xi\cdot\xi^t-M_{\gamma}^c\,.
\end{equation}
(We view $\xi$ as a column matrix.)
\end{clm}
\begin{proof}
By~\eqref{mucca} the ideal of $X_{\gamma}$  is generated by $\det M_{\gamma}$ and the entries of the matrices
in~\eqref{profumo}.
By Cramer's formula $\det M_{\gamma}$ belongs to the ideal generated by the entries of the two  matrices. This proves that the ideal of $X_{\gamma}$ is as claimed. 
\end{proof}
Now we suppose in addition that $R$ is a finitely generated $\CC$-algebra.
Let $p\in \Spec R$ be a closed point:   we are interested in the localization of $X_{\gamma}$ at points in $f_{\gamma}^{-1}(p)$.
Let $J\subset U^{\vee}(p)$ be a  subspace complementary to $\ker\gamma(p)$. Let ${\bf J}\subset U^{\vee}$ be a free submodule whose fiber over  $p$ is equal to $J$. Let ${\bf K}\subset U^{\vee}$ be the submodule orthogonal to ${\bf J}$ i.e.
\begin{equation}
{\bf K}:=\{u\in U^{\vee}\mid \gamma(a)(u)=0\quad\forall a\in {\bf J}\}\,.
\end{equation}
The localization of ${\bf K}$ at $p$  is  free. Let  
 $K:={\bf K}(p)$ be the fiber of ${\bf K}$ at $p$; clearly $K=\ker\gamma(p)$.  Localizing at $p$ we have 
 \begin{equation}\label{ugeikappa}
U^{\vee}_p={\bf K}_p\oplus {\bf J}_p\,. 
\end{equation}
Corresponding to~\eqref{ugeikappa} we may write 
 $ \gamma_p=\gamma_{\bf K}\oplus_{\bot} \gamma_{\bf J}$
 where $\gamma_{\bf K}\colon {\bf K}_p\to {\bf K}^{\vee}_p$ and $\gamma_J\colon {\bf J}_p\to {\bf J}^{\vee}_p$   are symmetric maps. Notice that we have an equality of germs 
 \begin{equation}\label{yuguali}
 (Y_{\gamma},p)= (Y_{\gamma_{\bf K}},p).
\end{equation}
We claim that there is a compatible isomorphism of germs $(X_{\gamma_{\bf K}},f_{\gamma_{\bf K}}^{-1}(p))\cong
(X_{\gamma},f_{\gamma}^{-1}(p))$.  In fact let $k:=\dim K$ and $d:=\rk\,U$. Choose bases of ${\bf K}_p$ and ${\bf J}_p$; by~\eqref{ugeikappa} we get a basis of $U^{\vee}_p$. The dual bases of ${\bf K}^{\vee}_p$, ${\bf J}^{\vee}_p$ and $U^{\vee}_p$ are compatible with respect to the decomposition dual to~\eqref{ugeikappa}. Corresponding to the chosen bases we have embeddings
 $X_{\gamma_K}\hra Y_{\gamma_K}\times\CC^k$ and 
$X_{\gamma}\hra Y_{\gamma}\times\CC^d$.
 The decomposition dual to~\eqref{ugeikappa} gives an embedding 
 $ j\colon Y_{\gamma_K}\times\CC^k\hra
Y_{\gamma}\times\CC^d$. 
 \begin{clm}\label{clm:nucleo}
Keep notation as above. The composition
\begin{equation}
X_{\gamma_K}\hra (Y_{\gamma_K}\times\CC^k)\overset{j}{\lra}
(Y_{\gamma}\times\CC^d)
\end{equation}
defines an isomorphism of  germs in the analytic topology
\begin{equation}
(X_{\gamma_{\bf K}},f_{\gamma_{\bf K}}^{-1}(p))\overset{\sim}{\lra}
(X_{\gamma},f_{\gamma}^{-1}(p))
\end{equation}
which commutes with the maps $f_{\gamma_{\bf K}}$ and $f_{\gamma}$.
\end{clm}
\begin{proof}
This follows by writing 
 $ \gamma_p=\gamma_{\bf K}\oplus_{\bot} \gamma_{\bf J}$ and by recalling~\eqref{mucca}. We pass to the analytic topology in order to be able to extract the square root of a regular non-zero function.
\end{proof}
\begin{prp}\label{prp:gradodue}
Assume that $R$ is a finitely generated $\CC$-algebra. 
Suppose that we have Exact Sequence~\eqref{mappasimm}.
 Then the following hold:
\begin{itemize}
\item[(1)]
$f_{\gamma}^{-1}Y_{\gamma}(1)\to Y_{\gamma}(1)$ is a topological covering of degree $2$.
\item[(2)]
Let $p\in (Y_{\gamma}\setminus Y_{\gamma}(1))$ be a closed point. The fiber $f_{\gamma}^{-1}(p)$ consists of a single point $q$.
Let $\xi_i$ be the coordinates on $X_{\gamma}$ associated to Embedding~\eqref{esplicito}; then  $\xi_i(q)=0$ for $i=1,\ldots,d$. 
\end{itemize}
\end{prp}
\begin{proof}
(1):  Localizing at $p\in Y_{\gamma}(1)$ and applying~\Ref{clm}{nucleo} we get Item~(1).  (2):  Since $\cork M_{\gamma}(p)\ge 2$ we have $M_{\gamma}^c(p)=0$. Thus  Item~(2) follows from~\Ref{clm}{ideawu}.
\end{proof}
We may associate a double cover   $f_{\gamma}\colon X_{\gamma}\to Y_{\gamma}$ to a map $\beta$ which is symmetric in the derived category.
\begin{hyp}\label{hyp:ipotizzo}
We have~\eqref{unozeroenne} with  $\alpha$  an isomorphism  and in addition $\alpha=\mu$. 
\end{hyp}
\begin{prp}\label{prp:redux}
Assume that~\Ref{hyp}{ipotizzo} holds. Then $R/(\det M_{\lambda})\oplus N$ equipped with the product given by~\eqref{moltiplico}  is a commutative (associative) ring.
\end{prp} 
\begin{proof}
Let
$\gamma:=\lambda\circ\mu^{-1}$ and $U:=U_0$. Then~\eqref{mappasimm} holds and 
the product defined by $m_{\beta}$ is equal to the product defined by $m_{\gamma}$. By~\Ref{prp}{associo}  we get  that  $R/(\det M_{\lambda})\oplus N$ is a commutative associative ring. 
\end{proof}
\begin{dfn}\label{dfn:simmetrizzo}
Suppose that~\Ref{hyp}{ipotizzo} holds: the {\it symmetrization of~\eqref{unozeroenne}} is Exact Sequence~\eqref{mappasimm} with $\gamma$ and $U$ as in the proof of~\Ref{prp}{redux}. 
\end{dfn}
\subsection{Structure sheaf of double EPW-sextics}\label{subsec:voiladouble}
\setcounter{equation}{0}
Let $A\in\lagr$ and suppose that $Y_A\not=\PP(V)$. We will define the associated double cover $X_A\to Y_A$ by applying the results of~\Ref{subsec}{prodenne}. 
Since $A$ is Lagrangian the symplectic form defines a canonical isomorphism $\bigwedge^3 V/A\cong A^{\vee}$; thus~\eqref{diecidieci} defines  a map of vector-bundles $\lambda_A\colon F\to A^{\vee}\otimes\cO_{\PP(V)}$.
Let $i\colon Y_A\hra\PP(V)$ be the inclusion map: 
since a local generator of $\det\lambda_A$ annihilates $\coker  (\lambda_A)$  there is a unique sheaf $\zeta_A$ on $Y_A$ such that  we have an exact sequence
\begin{equation}\label{eccozeta}
0\lra F\overset{\lambda_A}{\lra} A^{\vee}\otimes\cO_{\PP(V)}\lra
i_{*}\zeta_A\lra 0.
\end{equation}
Choose $B\in\lagr$ transversal to $A$. Thus we have a direct-sum decomposition
$\bigwedge^3 V=A\oplus B$  and hence a projection map  $\bigwedge^3 V\to A$ inducing   a map $ \mu_{A,B}\colon F\to A\otimes\cO_{\PP(V)}$. 
 We claim that there is a  commutative diagram with exact rows
\begin{equation}\label{spqr}
\begin{array}{ccccccccc}
0 & \to & F&\mapor{\lambda_A}& A^{\vee}\otimes\cO_{\PP(V)} & \lra & i_{*}\zeta_A
&
\to & 0\\
 & & \mapver{\mu_{A,B}}& &\mapver{\mu^{t}_{A,B}} &
&
\mapver{\beta_{A}}& & \\
0 & \to & A\otimes\cO_{\PP(V)}& \mapor{\lambda_A^{t}}& F^{\vee} & \lra &
Ext^1(i_{*}\zeta_A,\cO_{\PP(V)}) & \to & 0\,.
\end{array}
\end{equation}
In fact the second row is obtained by applying the $Hom(\,\cdot\, ,\cO_{\PP(V)})$-functor to~(\ref{eccozeta}) and the equality $\mu_{A,B}^t\circ\lambda_A=\lambda_A^t\circ\mu_{A,B}$ holds because
 $F$ is a Lagrangian sub-bundle of $\bigwedge^3 V\otimes\cO_{\PP(V)}$.
 Lastly $\beta_A$ is defined to be the unique map making the diagram commutative; it exists because the rows are exact. Notice that the map $\beta_A$ is independent of the choice of $B$ as  suggested by the notation. 
Next by applying the $Hom(i_{*}\zeta_A,\,\cdot\,)$-functor to the exact sequence
\begin{equation}
0\lra\cO_{\PP(V)}\lra\cO_{\PP(V)}(6)\lra\cO_{Y_A}(6)\lra 0
\end{equation}
 we get the exact sequence
\begin{equation}\label{romaleone}
0\lra i_{*}Hom(\zeta_A,\cO_{Y_A}(6))\overset{\partial}{\lra}
Ext^1(i_{*}\zeta_A,\cO_{\PP(V)})\overset{n}{\lra} Ext^1(i_{*}\zeta_A,\cO_{\PP(V)}(6))
\end{equation}
where $n$ is locally equal to multiplication  by $\det\lambda_A$. Since the second row of~(\ref{spqr}) is exact  a local generator of $\det\lambda_A$ annihilates  $Ext^1(i_{*}\zeta_A,\cO_{\PP(V)})$; thus $n=0$ and hence we get a canonical isomorphism
\begin{equation}\label{extugualehom}
\partial^{-1}\colon 
Ext^1(i_{*}\zeta_A,\cO_{\PP(V)})\overset{\sim}{\lra} i_{*}Hom(\zeta_A,\cO_{Y_A}(6)).
\end{equation}
We define  $\wt{m}_A$ by setting
\begin{equation}\label{emmea}
\begin{matrix}
\zeta_A\times\zeta_A & \overset{\wt{m}_A}\lra & \cO_{Y_A}(6) \\
(\sigma_1,\sigma_2) & \mapsto & (\partial^{-1}\circ\beta_A(\sigma_1))(\sigma_2).
\end{matrix}
\end{equation}
Let $\xi_A:=\zeta_A(-3)$. Tensorizing both sides of~\eqref{emmea} by $\cO_{Y_A}(-6)$ we get   a multiplication map 
\begin{equation}\label{aridaje}
\xi_A\times\xi_A  \overset{m_A}\lra  \cO_{Y_A}.
\end{equation}
 Thus we have defined a multiplication map on  $\cO_{Y_A}\oplus\xi_A$. The following result is well-known to experts. 
\begin{prp}\label{prp:mappabeta}
Let $A\in\lagr$ and suppose that $Y_A\not=\PP(V)$.  
Let notation be as above. Then:
\begin{itemize}
\item[(1)]
$\beta_A$ is an isomorphism.
\item[(2)]
The multiplication map $m_A$ is associative and commutative.
\end{itemize}
\end{prp}
\begin{proof}
Let $[v_0]\in\PP(V)$. Choose $B\in\lagr$ transversal to $F_{v_0}$ (and to $A$ of course). Then $\mu_{A,B}$ is an isomorphism in an open neighborhood $U$ of $[v_0]$. It follows that $\beta_A$ is an isomorphism in a neighborhood of $[v_0]$. This proves Item~(1). Let's prove Item~(2). Let $B\in\lagr$ and $U$ be as above; we may assume that $U$ is affine. Let $N:=H^0(i_{*}\zeta_A|_U)$ and 
$\beta:=H^0(\beta_A|_U)$. Thus $\beta\colon N\to {\rm Ext}^1(N,\CC[U])$.  
By Commutativity of Diagram~\eqref{spqr} and by~\Ref{prp}{redux} we get that the multiplication map $m_{\beta}$ is associative and commutative. On the other hand $m_{\beta}$ is the multiplication induced by $m_A$ on $N$; since $[v_0]$ is an arbitrary point of $\PP(V)$ it follows that  $m_A$ is associative and commutative. 
\end{proof}
We let $X_A:=\Spec(\cO_{Y_A}\oplus\xi_A)$ and  we let
$f_A\colon X_A\to Y_A$
be the structure morphism. Then $X_A$ is the {\it double EPW-sextic} associated to $A$ and $f_A$ is its structure map. 
The {\it covering involution} of $X_A$ is the automorphism $\phi_A\colon X_A\to X_A$ corresponding to the involution of   
$\cO_{Y_A}\oplus\xi_A$ with $(-1)$-eigensheaf equal to  $\xi_A$.
\subsection{Local models of double covers}\label{subsec:doppipiccoli}
\setcounter{equation}{0}
In the present subsection we assume that $R$ is a finitely generated  $\CC$-algebra. Let $\cW$ be a finite-dimensional complex vector-space.  We will suppose that  we have an exact sequence
\begin{equation}\label{ancoramu}
0\lra R\otimes\cW^{\vee}\overset{\gamma}{\lra} R\otimes\cW\lra N\lra 0,\qquad 
\gamma=\gamma^t.
\end{equation}
Thus we have a double cover $f_{\gamma}\colon X_{\gamma}\to Y_{\gamma}$.
Let $p\in Y_{\gamma}$ be a closed point. 
We will examine $X_{\gamma}$ in a neighborhood of $f_{\gamma}^{-1}(p)$ when the corank of $\gamma(p)$ is small. 
We may view $\gamma$ as a regular map $\Spec R\to \Sym^2 \cW$; thus it makes sense to consider the differential 
\begin{equation}
d\gamma(p)\colon T_p \Spec R\to \Sym^2\cW.
\end{equation}
 Let  $K(p):=\ker\gamma(p)\subset\cW^{\vee}$; we will consider the linear map
 \begin{equation}\label{difftre}
\begin{matrix}
T_p \Spec R & \overset{\delta_{\gamma}(p)}{\lra} &
\Sym^2 K(p)^{\vee} \\
\tau & \mapsto & d\gamma(p)(\tau)|_{K(p)}\,.
\end{matrix}
\end{equation}
 Let $d:=\dim\cW$; choosing a basis of $\cW$ 
we realize $X_{\gamma}$ as a subscheme of $\Spec R\times\CC^d$ with ideal given by~\Ref{clm}{ideawu}. Since $\cork\gamma(p)\ge 2$  \Ref{prp}{gradodue} gives that $f_{\gamma}^{-1}(p)$ consists of  a single point $q$ - in fact the $\xi_i$-coordinates of $q$ are all zero. Throughout this subsection we let
\begin{equation}
f_{\gamma}^{-1}(p)=\{q\}\,.
\end{equation}
\begin{clm}\label{clm:pocheq}
Keep notation  as above. Suppose  that $d=\dim\cW=2$ and that $\gamma(p)=0$. Then $I(X_{\gamma})$ is generated by the entries of $\xi\cdot\xi^t-M^c_{\gamma}$. 
\end{clm}
\begin{proof}
\Ref{clm}{ideawu} together with  a straightforward computation. 
\end{proof}
\begin{expl}\label{expl:rangodue}
{\rm Let $R=\CC[x,y,z]$, $\cW=\CC^2$. Suppose that the matrix associated to $\gamma$ is 
\begin{equation}
M_{\gamma}=
\left(
\begin{matrix}
x & y \\
y & z
\end{matrix}
\right).
\end{equation}
Then $f_{\gamma}\colon X_{\gamma}\to Y_{\gamma}$ is identified with
\begin{equation}
\begin{matrix}
\CC^2 & \lra & V(xz-y^2) \\
(\xi_1,\xi_2) & \mapsto & (\xi_2^2,\,-\xi_1\xi_2,\,\xi_1^2)
\end{matrix}
\end{equation}
i.e.~the quotient map for the action of $\la -1\ra$ on $\CC^2$.}
\end{expl}
\begin{prp}\label{prp:biancabernie}
Keep notation  as above. Suppose that 
the following hold:
\begin{itemize}
\item[(a)]
$\cork\,\gamma(p)=2$,
\item[(b)]
the localization $R_p$ is  regular.
\end{itemize}
Then $X_{\gamma}$ is smooth at $q$ if and only if   $\delta_{\gamma}(p)$ is surjective. 
\end{prp}
\begin{proof}
Applying~\Ref{clm}{nucleo} we get that we may assume that $d=2$. Let
\begin{equation}
M_{\gamma}=
\left(
\begin{matrix}
a & b \\
b & c
\end{matrix}
\right).
\end{equation}
By~\Ref{clm}{pocheq} the ideal of $X_{\gamma}$ in $\Spec R\times\CC^2$ is generated by the entries of $\xi\cdot\xi^t-M^c_{\gamma}$ i.e.
\begin{equation}\label{lassie}
I(X_{\gamma})=(\xi_1^2-c,\ \xi_1\xi_2+b,\ \xi_2^2-a)\,.
\end{equation}
Thus 
\begin{equation}
\cod(T_q X_{\gamma},T_q (\Spec R\times\CC^2))=\dim\la da(p)\,,db(p)\,,dc(p)\ra\,.
\end{equation}
On the other hand $\cod_q(X_{\gamma}\,,\Spec R\times\CC^2)=3$ and hence we get that $X_{\gamma}$ is smooth at $q$ if and only if $\delta_{\gamma}(p)$ is surjective. 
\end{proof}
\begin{clm}\label{clm:sempresing}
Keep notation and hypotheses as above. Suppose that $\cork\gamma(p)\ge 3$. Then $X_{\gamma}$ is  singular at $q$.
\end{clm}  
\begin{proof}
Let $I$ be the ideal of $X_{\gamma}$ in $\Spec R[\xi_1,\ldots,\xi_d]$.  By~\Ref{clm}{ideawu} we get that $I$   is non-trivial but the differential at $q$ of an arbitrary $g\in I$ is zero.
\end{proof}
Next we will discuss in greater detail those $X_{\gamma}$ whose corank at $f^{-1}_\gamma(p)$ is equal to $3$. 
First we will identify the \lq\lq universal\rq\rq example (the universal example for corank $2$ is~\Ref{expl}{rangodue}). Let $\cV$ be a $3$-dimensional complex vector space. We view $\Sym^2\cV$ as an affine ($6$-dimensional) space and we let $R:=\CC[\Sym^2\cV]$ be its ring of regular functions. 
We identify $R\otimes_{\CC}\cV$ and $R\otimes_{\CC}\cV^{\vee}$ with the space of $\cV$-valued, respectively  $\cV^{\vee}$-valued, regular maps on  $\Sym^2\cV$.
Let
\begin{equation}
R\otimes_{\CC}\cV^{\vee}  \overset{\gamma}{\lra} 
R\otimes_{\CC}\cV 
\end{equation}
be the map induced on the spaces of global sections by the tautological
map of vector-bundles $\Spec R\times\cV^{\vee}\lra \Spec R\times\cV$.   The map $\gamma$ is  symmetric. Let $N$ be the cokernel of $\gamma$: thus 
\begin{equation}
0\lra R\otimes_{\CC}\cV^{\vee} \overset{\gamma}{\lra} 
R\otimes_{\CC}\cV\lra N\lra 0
\end{equation}
is an exact sequence. Since $\gamma$ is  symmetric  it defines a double cover $f\colon X(\cV)\to Y(\cV)$ where
\begin{equation}
Y(\cV):=\{\alpha\in \Sym^2\cV\mid \rk\,\alpha<3 \}
\end{equation}
is the variety of degenerate quadratic forms. We let 
\begin{equation}\label{invoix}
\phi\colon X(\cV)\to X(\cV)
\end{equation}
 be the covering involution of $f$. 
One describes explicitly $X(\cV)$ as follows. Let  
\begin{equation}
(\cV\otimes\cV)_1:= \{\mu\in (\cV\otimes\cV)\mid \rk\,\mu \le 1\}.
\end{equation}
Thus $(\cV\otimes\cV)_1$ is the cone over the Segre variety $\PP(\cV)\times\PP(\cV)$. We have a finite degree-$2$ map 
\begin{equation}\label{mappasigma}
\begin{matrix}
(\cV\otimes\cV)_1 & \overset{\sigma}{\lra} & Y(\cV) \\
\mu & \mapsto & \mu+\mu^t\,.
\end{matrix}
\end{equation}
\begin{prp}\label{prp:equarela}
Keep notation as above. There exists  a commutative diagram 
\begin{equation}\label{triangolo}
\xymatrix{  \\  (\cV\otimes\cV)_1\ar^{\sigma}[dr]  \ar^{\tau}[rr] & & X(\cV)\ar_{f}[dl] \\
& Y(\cV) &}
\end{equation}
where $\tau$ is an isomorphism. Let $\phi$ be Involution~\eqref{invoix}: then 
\begin{equation}\label{invtras}
\phi\circ\tau(\mu)=\tau(\mu^t),\qquad \forall\mu\in (\cV\otimes\cV)_1.
\end{equation}
\end{prp}
\begin{proof}
In order to define $\tau$ we will give a coordinate-free version of Inclusion~(\ref{esplicito}) in the case of $X(\cV)$.
Let
\begin{equation}\label{spigola}
\begin{matrix}
\Sym^2\cV\times(\cV^{\vee}\otimes\bigwedge^3\cV) &
\overset{\Psi}{\lra} & (\cV\otimes\bigwedge^3\cV)\times
(\cV^{\vee}\otimes\cV^{\vee}\otimes\bigwedge^3\cV\otimes\bigwedge^3\cV) 
\\
(\alpha,\xi) & \mapsto & (\alpha\circ\xi,\ \xi^t\circ\xi - \bigwedge^ 2 \alpha)\,.
\end{matrix}
\end{equation}
A few words of explanation. 
In the definition of the first component of $\Psi(\alpha,\xi)$ we view  $\xi$ as belonging to $\Hom(\bigwedge^3\cV^{\vee},\cV^{\vee})$, in   the definition of the second component of $\Psi(\alpha,\xi)$ we view  $\xi$ as belonging to 
$ \Hom(\cV\otimes\bigwedge^3\cV^{\vee},\CC)$. Moreover we make the obvious choice of isomorhpism $\CC\cong\CC^{\vee}$.  Secondly 
\begin{equation}
\bigwedge^ 2\alpha\in\Hom(\bigwedge^2\cV^{\vee},\,\bigwedge^2\cV)
=\Hom(\cV\otimes\bigwedge^3\cV^{\vee},\,
\cV^{\vee}\otimes\bigwedge^3\cV)=
\cV^{\vee}\otimes\cV^{\vee}\otimes
\bigwedge^3\cV\otimes\cV\,.
\end{equation}
Choosing a basis of $\cV$ we get an embedding $X(\cV)\subset \Sym^2\cV\times \CC^3$, see~\eqref{esplicito}. \Ref{clm}{ideawu} gives equality of pairs
\begin{equation}\label{coppie}
(\Sym^2\cV\times(\cV^{\vee}\otimes\bigwedge^3\cV),\,
\Psi^{-1}(0))= (\Sym^2\cV\times\CC^3,\,X(\cV))\,,
\end{equation}
where $\Psi^{-1}(0)$ is the scheme-theoretic fiber of $\Psi$. 
Now notice that we have an isomorphism
\begin{equation}\label{cuffia}
\begin{matrix}
\cV\otimes\cV & \overset{\overset{T}{\sim}}{\lra} &
 \Sym^2\cV\times(\cV^{\vee}\otimes\bigwedge^3\cV) \\
  \epsilon & \mapsto & (\epsilon+\epsilon^t,\epsilon-\epsilon^t)\,.
\end{matrix}
\end{equation}
Let $\tau:=T|_{(\cV\otimes\cV)_1}$: thus we have an embedding 
\begin{equation}\label{eccotau}
\tau\colon(\cV\otimes\cV)_1\hra \Sym^2\cV\times(\cV^{\vee}\otimes\bigwedge^3\cV).
\end{equation}
  We will show that we have
equality of  schemes
\begin{equation}\label{ixgielle}
\im(\tau)=\Psi^{-1}(0)(=X(\cV)).
\end{equation}
 First let
\begin{equation}
\begin{matrix}
\cV\oplus\cV &
\overset{\rho}{\lra} & (\cV\otimes\cV)_1  \\
(\eta,\beta) & \mapsto & \eta^t\circ\beta.
\end{matrix}
\end{equation}
Notice that $\rho$ is the quotient map for the $\CC^{\times}$-action on $\cV\oplus\cV$ defined by $t(\eta,\beta):=(t\eta,t^{-1}\beta)$.  We  have
\begin{equation}
\tau\circ\pi= (\eta^t\circ\beta+\beta^t\circ\eta, 
 \eta\wedge\beta).
\end{equation}
Let's prove that 
\begin{equation}\label{tauinix}
\Psi^{-1}(0)\supset\im(\tau)\,.
\end{equation}
Notice that $\GL(\cV)$ acts on  $(\cV\otimes\cV)_1$ with a unique dense orbit namely $\{\eta^t\circ \beta\mid \eta\wedge\beta\not=0\}$. An easy computation shows that $\tau(\eta^t\circ \beta)\in \Psi^{-1}(0)$ for  a conveniently chosen $\eta^t\circ \beta$ in the dense orbit of $(\cV\otimes\cV)_1$;  it follows that~\eqref{tauinix} holds. On the other hand $T$ defines an isomorphism of pairs
\begin{equation}\label{coppiedue}
(\cV\otimes\cV,\, (\cV\otimes\cV)_1)\cong
 (\Sym^2\cV^{\vee}\times (\cV^{\vee}\otimes\bigwedge^3\cV),\,\im(\tau))\,.
\end{equation}
Since the ideal of $(\cV\otimes\cV)_1$ in $\cV\otimes\cV$  is generated by $9$ linearly independent quadrics we get that the ideal of $\im(\tau)$ in $\Sym^2\cV^{\vee}\times (\cV^{\vee}\otimes\bigwedge^3\cV)$ is generated by $9$ linearly independent quadrics. The ideal of $\Psi^{-1}(0)$ in $\Sym^2\cV\times(\cV^{\vee}\otimes\bigwedge^3\cV)$ is likewise generated by $9$ linearly independent quadrics  - see~\eqref{spigola}. Since 
$\Psi^{-1}(0)\supset\im(\tau)$ we get that the ideals of $\Psi^{-1}(0)$ and of $\im(\tau)$ are the same and hence~\eqref{ixgielle} holds.
 This proves that $\tau$ is an isomorphism between $(\cV\otimes\cV)_1$ and $X(\cV)$. Diagram~\eqref{triangolo} is commutative by construction. 
 Equation~\eqref{invtras}  is equivalent to the equality 
\begin{equation}
\phi(\tau\circ\rho(\beta,\eta))=\tau\circ \rho(\eta,\beta)).
\end{equation}
The above equality holds because $\beta\wedge\eta=-\eta\wedge\beta$.
\end{proof}
The following result  is an immediate consequence of~\Ref{prp}{equarela}.
\begin{crl}\label{crl:singixvu}
$\sing X(\cV)=\tau(0)=f^{-1}(0)$.
\end{crl}
\section{The divisor $\Delta$}\label{sec:avoidelta}
\setcounter{equation}{0}
\subsection{Parameter counts}
\setcounter{equation}{0}
Let $\Delta_{+}\subset\lagr$ and $\wt{\Delta}_{+},\wt{\Delta}_{+}(0)\subset\lagr\times\PP(V)^2$ be  
\begin{eqnarray}
\Delta_{+}:= & \{A\in\lagr\mid | Y_A[3]|>1\}\,,
\index{$\Delta_{+}$}\\
\wt{\Delta}_{+}:=& \{(A,[v_1],[v_2])\mid [v_1]\not=[v_2],
\quad \dim(A\cap F_{v_i})\ge 3\}\,,\\
\wt{\Delta}_{+}(0):=& \{(A,[v_1],[v_2])\mid [v_1]\not=[v_2],
\quad \dim(A\cap F_{v_i})= 3\}\,.
\end{eqnarray}
Notice that $\wt{\Delta}_{+}$ and $\wt{\Delta}_{+}(0)$ are locally closed. 
\begin{lmm}\label{lmm:codelta}
Keep notation as above. The following hold:
\begin{enumerate}
\item[(1)]
$\wt{\Delta}_{+}$ is irreducible of dimension $53$. 
\item[(2)]
 $\Delta_{+}$ is   constructible and $\cod(\Delta_{+},\lagr)\ge 2$. 
\end{enumerate}
\end{lmm}
\begin{proof}
(1): Let's prove that $\wt{\Delta}_{+}(0)$ is irreducible of dimension $53$.  
Consider the map
\begin{equation}
\begin{matrix}
\wt{\Delta}_{+}(0) & \overset{\eta}{\lra} & 
\Gr(3,\bigwedge^3 V)^2\times\PP(V)^2 \\
(A,[v_1],[v_2]) & \mapsto & 
(A\cap F_{v_1},A\cap F_{v_2},[v_1],[v_2])\,.
\end{matrix}
\end{equation}
We have
\begin{equation}
\im\eta=\{(K_1,K_2,[v_1],[v_2])\mid 
 K_i\in\Gr(3,F_{v_i}),\quad K_1\bot K_2,\quad [v_1]\not=[v_2]\}\,.
\end{equation}
We stratify $\im\eta$  according to $i:=\dim(K_1\cap F_{v_2})$ and to $j:=\dim(K_1\cap K_2)$; of course $j\le i$. Let $(\im\eta)_{i,j}\subset \im\eta$ be the stratum corresponding to $i,j$. A straightforward computation gives that
\begin{multline}
\dim\eta^{-1}(\im\eta)_{i,j}=10+7(3-i)+j(i-j)+(3-j)(4+i)+
\frac{1}{2}(j+5)(j+4)=\\
=53-4i-\frac{1}{2}j(j-1)\,.
\end{multline}
Since $0\le i, j$ one gets that the maximum is achieved for $i=j=0$ and that it equals $53$. It follows that $\wt{\Delta}_{+}(0)$ is irreducible of dimension $53$. On the other hand $\wt{\Delta}_{+}(0)$ is dense in $\wt{\Delta}_{+}$ (easy) and hence we get that Item~(1) holds. 
  (2): Let $\pi_{+}\colon \wt{\Delta}_{+}\to \lagr$ be the forgetful map: $\pi_{+}([v_1],[v_2],A)=A$. Then $\pi_{+} (\wt{\Delta}_{+})=\Delta_{+}$. By Item~(1) we get that $\dim\Delta_{+}\le 53$: since $\dim\lagr=55$ we get that Item~(2) holds.
\end{proof}
\begin{prp}\label{prp:codelta}
 The following hold:
\begin{itemize}
\item[(1)]
 $\Delta$ is  closed irreducible of codimension $1$ in $\lagr$ and not equal to $\Sigma$.  
\item[(2)]
If $A\in\Delta$ is generic then $Y_A[3]=Y_A(3)$ and it consists of a single point.
\end{itemize}
\end{prp}
\begin{proof}
(1): Let 
\begin{equation}
\wt{\Delta}:=\{(A,[v])\mid \dim(F_v\cap A)\ge 3\},\qquad 
\wt{\Delta}(0):=\{(A,[v])\mid \dim(F_v\cap A)= 3\}.
\end{equation}
Then $\wt{\Delta}$ is a closed subset of  $\lagr\times\PP(V)$ and $\wt{\Delta}(0)$ is an open subset of  $\wt{\Delta}$. Let 
$\pi\colon\wt{\Delta}\to\lagr$
 be the forgetful map. Thus 
$ \pi(\wt{\Delta})=\Delta$: since $\pi$ is projective it follows that $\Delta$ is closed. 
Projecting $\wt{\Delta}(0)$ to $\PP(V)$ we get that $\wt{\Delta}(0)$ is smooth irreducible of dimension $54$. 
A standard dimension count shows that $\wt{\Delta}(0)$ is open dense in $\wt{\Delta}$; thus $\wt{\Delta}$ is irreducible  of dimension $54$. It follows that $\Delta$ is irreducible. 
By~\Ref{lmm}{codelta} we know that $\dim\wt{\Delta}_{+}\le 53$.  
 It follows that the generic fiber of $\wt{\Delta}\to\Delta$ is a single point, in particular $\dim\Delta=54$ and hence $\cod(\Delta,\lagr)=1$ because $\dim\lagr=55$. A dimension count shows that $\dim(\Delta\cap\Sigma)<54$ and hence $\Delta\not=\Sigma$. This finishes the proof of Item~(1). (2): Let $A\in\Delta$ be generic: we already noticed that there exists a unique $[v]\in\PP(V)$ such that $([v],A)\in \wt{\Delta}$, i.e.~$Y_A[3]$ consists of a single point. Since $\wt{\Delta}(0)$ is dense in $\wt{\Delta}$ and $\dim\wt{\Delta}=\dim\Delta$ we get that  $([v],A)\in \wt{\Delta}(0)$, i.e.~$Y_A[3]=Y_A(3)$. This finishes the proof of Item~(2). 
\end{proof}
\subsection{First order computations}
\setcounter{equation}{0}
 Let  $(A,[v_0])\in\wt{\Delta}(0)$. We will study the differential of $\pi\colon\wt{\Delta}\to\lagr$
 at $(A,[v_0])$. 
 First we will give a  local description of $\wt{\Delta}$ as  degeneracy locus. Let  
\begin{equation}\label{ennevu}
\NN(V):=\{A\in\lagr \mid Y_A=\PP(V)\}.
\end{equation}
Notice that $\NN(V)$ is closed.  
Let $\cY$  be the tautological family of EPW-sextics
i.e.
\begin{equation}
\cY:=\{(A,[v])\in (\lagr\setminus\NN(V))\times\PP(V)\mid 
\dim(A\cap F_v)>0\}\,.
\end{equation}
Of course $\cY$ has a description as a determinantal variety 
 and hence it has a 
natural scheme structure. For $\cU\subset (\lagr\setminus\NN(V))$  open we let $\cY_{\cU}:=\cY\cap(\cU\times\PP(V))$.
Given $B\in \lagr$  let
 \begin{equation}\label{eccoua}
U_B:=\{A\in \lagr \mid A\pitchfork B\}\setminus\NN(V).
\end{equation}
(Here $A\pitchfork B$ means that $A$ intersects transversely $B$ i.e.~$A\cap B=\{0\}$.) Let $i_{U_B}\colon\cY_{U_B}\hra U_B\times\PP(V)$ be the inclusion and let $\cA$ be the tautological rank-$10$ vector-bundle on $\lagr$ (the fiber of $\cA$ over $A$ is $A$ itself).
Going through the argument that 
 produced Commutative Diagram~\eqref{spqr} we get that there exists a commutative diagram
\begin{equation}\label{pazzaidea}
\begin{array}{ccccccccc}
0 & \to & \cO_{U_B}\boxtimes F&\mapor{\lambda_{U_B}}&
(\cA^{\vee}|_{U_B})\boxtimes\cO_{\PP(V)} & 
\lra & i_{U_B,*}\zeta_{U_B} &
\to & 0\\
 & & \mapver{\mu_{U_B}}& &\mapver{\mu^{t}_{U_B}} &
&
\mapver{\beta_{U_B}}& & \\
0 & \to & (\cA|_{U_B})\boxtimes\cO_{\PP(V)}& 
\mapor{\lambda_{U_B}^{t}}&\cO_{U_B}\boxtimes F^{\vee} & \lra &
Ext^1(i_{U_B,*}\zeta_{U_B},\cO_{U_B\times\PP(V)}) & \to & 0
\end{array}
\end{equation}
Now let $(A,[v_0])\in\cY$. Choose $B\in\lagr$ such that $B\pitchfork A$ and $B\pitchfork F_{v_0}$.  Let $\cN\subset\PP(V)$ be an open neighborhood of $[v_0]$ such that $B\pitchfork F_w$ for all $w\in\cN$. The restriction to $U_B$ of  $\cA$    is trivial and  the restriction to $\cN$ of $F$ is likewise trivial.  Moreover the restriction of $\mu_{U_B}$ to $U_B\times\cN$ is an isomorphism. Let 
\begin{equation}\label{gamfam}
\gamma:=(\lambda_{U_B}|_{U_B\times\cN})\circ (\mu_{U_B}|_{U_B\times\cN})^{-1}. 
\end{equation}
We have an exact sequence
\begin{equation}\label{pattypravo}
0\lra (\cA|_{U_B})\boxtimes\cO_{\cN}\overset{\gamma}{\lra} (\cA^{\vee}|_{U_B})\boxtimes\cO_{\cN}
\lra i_{U_B,*}\zeta_{U_B} |_{U_B\times\cN} \lra 0
\end{equation}
  The map $\gamma$ is symmetric, in fact it is the symmetrization of the restriction of~\eqref{pazzaidea} to $U_B\times\cN$ - see~\Ref{dfn}{simmetrizzo}.  Then $\wt{\Delta}\cap(U_B\times \cN)$ is the symmetric degeneration locus
\begin{equation}\label{degsim}
\wt{\Delta}\cap(U_B\times \cN)=\{(A',[v])\in (U_B\times \cN)  \mid \cork\gamma(A',[v])\ge 3\}
\end{equation}
and hence it inherits a natural structure of closed subscheme of $\lagr\times\PP(V)$. In order to study the differential of the forgetful map $\wt{\Delta}\to\PP(V)$ we will introduce some notation.
 Given $v\in V$ we define 
a  quadratic form $\phi^{v_0}_v$ on $F_{v_0}$ as follows. Let $\alpha\in F_{v_0}$; then $\alpha=v_0\wedge\beta$ for some $\beta\in\bigwedge^2 V$. We set
\begin{equation}\label{quadricapluck}
\phi^{v_0}_v(\alpha):=\vol(v_0\wedge v\wedge\beta\wedge\beta).
\end{equation}
The above equation gives a well-defined quadratic form on $F_{v_0}$ because  $\beta$ is determined up to addition by an element of $F_{v_0}$. Of course $\phi^{v_0}_v$ depends only on the class of $v$ in $V/[v_0]$. Choose a direct-sum decomposition
\begin{equation}\label{scelgodeco}
V=[v_0]\oplus V_0.
\end{equation}
We have the isomorphism
\begin{equation}\label{trezeguet}
\begin{matrix}
\lambda^{v_0}_{V_0}\colon\bigwedge^2 V_0 & 
\overset{\sim}{\lra} & F_{v_0}\\
\beta & \mapsto & v_0\wedge\beta\,.
\end{matrix}
\end{equation}
Under the above identification  the Pl\"ucker quadratic forms on  $\bigwedge^2 V_0$ correspond to the quadratic forms $\phi^{v_0}_v$ for $v$ varying in $V_0$. 
 Let
 $K:=A\cap F_{v_0}$ and
\begin{equation}\label{botte}
\begin{matrix}
V_0 & \overset{\tau^{v_0}_K}{\lra} & \Sym^2 K^{\vee} \\
v & \mapsto & \phi^{v_0}_v|_K
\end{matrix}
\qquad\qquad
\begin{matrix}
\Sym^2 A^{\vee} & \overset{\theta^{A}_{K}}{\lra} & \Sym^2 K^{\vee} \\
q & \mapsto & q|_K\,.
\end{matrix}
\end{equation}
The isomorphism
\begin{equation*}
\begin{matrix}
V_0 & \overset{\sim}{\lra} & \PP(V)\setminus\PP(V_0) \\
v & \mapsto & [v_0+v]
\end{matrix}
\end{equation*}
defines an isomorphism 
$V_0\cong T_{[v_0]}\PP(V)$.
Recall that the tangent space to $\lagr$ at $A$ is canonically identified with  $\Sym^2 A^{\vee}$. 
\begin{prp}\label{prp:tandeltil}
Keep notation as above - in particular choose~\eqref{scelgodeco}. Then 
\begin{equation}
T_{(A,[v_0])}\wt{\Delta}\subset T_{(A,[v_0])}\left(\lagr\times\PP(V)\right)=
\Sym^2 A^{\vee}\oplus V_0
\end{equation}
 is given by
\begin{equation}
T_{([v_0],A)}\wt{\Delta}=
\{(q,v) \mid \theta^{A}_{K}(q)-\tau^{v_0}_K(v)=0\}.
\end{equation}
\end{prp}
\begin{proof}
By the (local) degeneracy description~\eqref{degsim} we get that $(q,v)\in T_{([v_0],A)}\wt{\Delta}$ if and only if 
\begin{equation*}
0=d\gamma(A,[v_0])(q,v)|_K=d\gamma(A,[v_0])(q,0)|_K+d\gamma(A,[v_0])(0,v)|_K.
\end{equation*}
It is clear that $d\gamma(A,[v_0])(q,0)|_K=\theta^{A}_{K}(q)$. On the other hand Equation~(2.26) of~\cite{og2} gives that
\begin{equation}\label{labambola}
d\gamma(A,[v_0])(0,v)|_K=-\tau^{v_0}_K(v). 
\end{equation}
 The proposition follows.
\end{proof}
\begin{crl}\label{crl:tanramdel}
$\wt{\Delta}(0)$ is smooth (of codimension $6$ in $\lagr\times\PP(V)$).  
Let $(A,[v_0])\in\wt{\Delta}(0)$ and  $K:=A\cap F_{v_0}$.  The
 differential $d\pi(A,[v_0])$ is injective if and only if $\tau^{v_0}_K$ is injective.
\end{crl}
\begin{proof}
Let $(A,[v_0])\in\wt{\Delta}(0)$ and $K:=A\cap F_{v_0}$. The map $\theta^A_K$ is surjective: by~\Ref{prp}{tandeltil} we get that $T_{(A,[v_0])}\wt{\Delta}(0)$ has codimension $6$ in $T_{(A,[v_0])}(\lagr\times\PP(V))$.
On the other hand the description of  $\wt{\Delta}(0)$ as a symmetric degeneration locus - see~\eqref{degsim} -  gives that $\wt{\Delta}(0)$ has codimension at most $6$ in $\lagr\times\PP(V)$: it follows that $\wt{\Delta}(0)$ is smooth of codimension $6$ in $\lagr\times\PP(V)$.  The statement about injectivity of  $d\pi(A,[v_0])$ follows at once form~\Ref{prp}{tandeltil}.
\end{proof}
A comment regarding~\Ref{crl}{tanramdel}. The statement  about smoothness of $\wt{\Delta}(0)$  is \emph{not} contained in the proof of~\Ref{prp}{codelta} because in that proof we consider  $\wt{\Delta}(0)$ with its reduced structure. Before stating the next result we give the following definition: given $A\in\lagr$ we let
\begin{equation}\label{eccoteta}
\Theta_A:=\{W\in\Gr(3,V) \mid \bigwedge^3 W\subset A\}.
\end{equation}
\begin{prp}\label{prp:aloha}
Let $(A,[v_0])\in\wt{\Delta}(0)$ and let $K:=A\cap F_{v_0}$. Then 
$\tau^{v_0}_K$ is injective if and only if 
\begin{itemize}
\item[(1)]
no $W\in\Theta_A$  contains $v_0$, or
\item[(2)]
there is exactly one $W\in\Theta_A$ containing $v_0$ and moreover 
\begin{equation}
A\cap F_{v_0}\cap (\bigwedge^2 W\wedge V)=\bigwedge^3 W.
\end{equation}
\end{itemize}
If Item~(1) holds then $\im\tau^{v_0}_K$ belongs to the unique open $\PGL(K)$-orbit of $\Gr(5, \Sym^2 K^{\vee})$,  if Item~(2) holds then  $\im\tau^{v_0}_K$  belongs to the unique closed $\PGL(K)$-orbit  of $\Gr(5, \Sym^2 K^{\vee})$.
\end{prp}
\begin{proof}
Let $V_0\subset V$ be a codimension-$1$ subspace transversal to $[v_0]$. Let 
\begin{equation}\label{rocco}
\rho^{v_0}_{V_0}\colon F_{v_0}\overset{\sim}{\lra} \bigwedge^2 V_0
\end{equation}
be the inverse of Isomorphism~\eqref{trezeguet}. 
Let ${\bf K}:=\PP(\rho^{v_0}_{V_0}(K))\subset\PP(\bigwedge^2 V_0)$; then ${\bf K}$ is a projective plane. Isomorphism $\rho^{v_0}_{V_0}$ identifies the space of quadratic forms $\phi^{v_0}_v$, for $v\in V_0$, with the space of Pl\"ucker quadratic forms on $\bigwedge^2 V_0$. Since the ideal of $\Gr(2,V_0)\subset\PP(\bigwedge^2 V_0)$ is generated by the Pl'ucker quadratic forms we get that $\tau^{v_0}_K$ is identified with the natural restriction map
\begin{equation}
V_0=H^0(\cI_{\Gr(2,V_0)}(2))\overset{\tau^{v_0}_K}{\lra} H^0(\cO_{{\bf K}}(2))=\Sym^2 K^{\vee}.
\end{equation}
It follows that if the  scheme-theoretic intersection ${\bf K}\cap \Gr(2,V_0)$ is not empty nor  a single reduced point then $\tau^{v_0}_K$ is not injective. Now suppose that ${\bf K}\cap \Gr(2,V_0)$ is 
\begin{itemize}
\item[($1'$)]
empty i.e.~Item~(1) holds, or
\item[($2'$)]
 a single reduced point, i.e.~Item~(2) holds. 
\end{itemize}
Let 
\begin{equation}
\PP(\bigwedge^2 V_0)\overset{\Phi}{\dashrightarrow}  | H^0(\cI_{\Gr(2,V_0)}(2))|^{\vee} =\PP(V_0^{\vee})
\end{equation}
be the natural map: it associates to $[\alpha]\notin\Gr(2,V_0)$ the projectivization of $\supp\alpha$. We have a tautological identification
\begin{equation*}
{\bf K}\overset{\Phi|_{\bf K}}{\dashrightarrow}\PP(\im \tau^{v_0}_K)^{\vee}
\end{equation*}
and of course $\Phi|_{\bf K}$ is the Veronese embedding ${\bf K}\to |\cO_{\bf K}(2)|^{\vee}$ followed by 
the projection with center $\PP(\Ann (\im \tau^{v_0}_K))$. Notice that $\tau^{v_0}_K$ is not injective if and only if $\dim\PP(\Ann (\im \tau^{v_0}_K))\ge 1$. Now suppose that~($1'$)  holds. Then $\Phi|_{\bf K}$ is regular and in fact it is an isomorphism onto its image - see Lemma~2.7 of~\cite{og5}. Since the chordal variety of the Veronese surface in $|\cO_{\bf K}(2)|^{\vee}$ is a hypersurface it follows that $\dim\PP(\Ann (\im \tau^{v_0}_K))<1$ and hence $\tau^{v_0}_K$ is injective.We also get that $\Ann (\im \tau^{v_0}_K)$ is a point in $|\cO_{\bf K}(2)|^{\vee}$ which does not belong to  the chordal variety of the Veronese surface and hence it 
 belongs to unique open $\PGL(K)$-orbit. 
 Now suppose that~($2'$) holds. 	Assume that  $\tau^{v_0}_K$ is not injective. Then $\dim\PP(\Ann (\im \tau^{v_0}_K))\ge 1$. It follows that there exist $[x]\not=[y]\in{\bf K}$ in the regular locus of $\Phi|_{\bf K}$ (i.e.~neither $x$ nor $y$ is decomposable)  such that $\Phi([x])=\Phi([y])$.   By the description of $\Phi$ given above in terms of supports we get that $\supp(x)=\supp(y)=U$ where $\dim U=4$; since $\Gr(2,U)$ is a hypersurface in $\PP(\bigwedge^2 U)$ we get that the line $\la [x],[y]\ra\subset\PP(\bigwedge^2 V_0)$ intersects $\Gr(2,U)$ in a subscheme of length $2$. Since $\la [x],[y]\ra\subset{\bf K}$ it follows that ${\bf K}\cap \Gr(2,V_0)$ contains a scheme of length $2$,  that contradicts Item~($2'$).  This  proves that if~($2'$) holds then $\tau^{v_0}_K$ is injective. It also follows that $\Ann(\tau^{v_0}_K)$ belongs to the Veronese surface in $ |\cO_{\bf K}(2)|^{\vee}$ i.e.~$\im(\tau^{v_0}_K)$ belongs to the unique closed  $\PGL(K)$-orbit. 
\end{proof}
\section{Simultaneous resolution}\label{sec:famdesing}
\setcounter{equation}{0}
In the first subsection we will analyze families of double EPW-sextics and their singular locus. The second subsection shows how to construct the simultaneous desingularization described in Item~(3) of~\Ref{sec}{prologo} (the relation with the Hilbert square of a $K3$ will be given in~\Ref{sec}{zitelle}).
\subsection{Families of double EPW-sextics}\label{subsec:dopfam}
\setcounter{equation}{0}
Let $\cU\subset(\lagr\setminus\NN(V))$ (see~\eqref{ennevu}) be open.  Suppose that there exist a scheme $\cX_{\cU}$ and  a finite $f_{\cU}\colon \cX_{\cU}\to \cY_{\cU}$ 
such that for every $A\in\cU$ the induced map $f^{-1}Y_A\to Y_A$ is identified with $f_A\colon X_A\to Y_A$:  
 then we say that a {\it tautological family of double EPW-sextics parametrized by $\cU$ exists} - often we simply state that
 $f_{\cU}\colon \cX_{\cU}\to \cY_{\cU}$  exists. Composing $f_{\cU}$ with the natural map $\cY_{\cU}\to\cU$ we get a map 
$\rho_{\cU}\colon \cX_{\cU}\to\cU$  such that $\rho_{\cU}^{-1}(A)\cong X_A$.  
\begin{prp}\label{prp:eccotaut}
Let  $B\in\lagr$. A tautological family of double EPW-sextics parametrized by $U_B$ exists ($U_B$ is given by~\eqref{eccoua}).
\end{prp}
\begin{proof}
Let  $\nu\colon \cY_{U_B}\to\PP(V)$ be projection.
Let $\xi_{U_B}:=\zeta_{U_B}\otimes\nu^{*}\cO_{\PP(V)}(-3)$ where $\zeta_{U_B}$ is the sheaf on $\cY_{U_B}$ fitting in~\eqref{pazzaidea}. Look at Commutative Diagram~\eqref{pazzaidea}: proceeding as in the definition of the multiplication on $\cO_{Y_A}\oplus\xi_A$ we get that $\beta_{U_B}$ defines a multiplication on $\cO_{\cY_{U_B}}\oplus\xi_{U_B}$. By~\Ref{prp}{redux} we get that $\cO_{\cY_{U_B}}\oplus\xi_{U_B}$ is an associative commutative ring.
Let $\cX_{U_B}:=\Spec(\cO_{\cY_{U_B}}\oplus\xi_{U_B})$ and  $f_{U_B}\colon\cX_{U_B}\to\cY_{U_B}$ be the structure map.  
\end{proof}
Let $\cU\subset (\lagr\setminus\NN(V))$ be open and such that $f_{\cU}\colon\cX_{\cU}\to \cY_{\cU}$ exists.
We will determine the singular locus of $\cX_{\cU}$. 
Let
\begin{eqnarray}
\cY[d]:= & \{(A,[v])\in(\lagr\setminus\NN(V))\times\PP(V)
  \mid \dim(A\cap F_v)\ge d\},\\
 \cY(d):= & \{(A,[v])\in(\lagr\setminus\NN(V))\times\PP(V)
  \mid \dim(A\cap F_v)= d\}. 
\end{eqnarray}
Then $\cY[d]$ has a natural structure of closed subscheme of $\lagr\times\PP(V)$ given by its local description as a symmetric determinantal variety - see {\bf Subsection~2.2} of~\cite{og5}. 
Let $\cU\in (\lagr\setminus\NN(V))$ be open. We let $\cY_{\cU}[d]:=\cY[d]\cap \cY_{\cU}$ and similarly for $\cY_{\cU}(d)$. Suppose that $f_{\cU}\colon\cX_{\cU}\to\cY_{\cU}$ is defined;  we let  
\begin{equation}
\cW_{\cU}:= f_{\cU}^{-1}\cY[3].
\end{equation}
Notice that the restriction of $f_{\cU}$ to  $\cW_{\cU}$ defines an isomorphism
$\cW_{\cU}\overset{\sim}{\lra}\cY_{\cU}[3]$. We will prove the following result.
\begin{prp}\label{prp:singfam}
Let $\cU\subset (\lagr\setminus\NN(V))$ be open and suppose that  $f_{\cU}\colon\cX_{\cU}\to\cY_{\cU}$ exists. Then $\sing\cX_{\cU}=\cW_{\cU}$.
\end{prp}
\begin{proof}
We may assume that $\cU=U_B\times\cN$ where $B\in\lagr$ and  $\cN\subset\PP(V)$ is an open subset such that $B\pitchfork F_w$ for all $w\in\cN$.  Then (see the proof of~\Ref{prp}{eccotaut}) 
\begin{equation}\label{pavia}
\text{$f_{U_B}^{-1}(\cU)=X_{\gamma}$ where $\gamma$ is given by~\eqref{gamfam}.}
\end{equation}
Thus it suffices to examine $X_{\gamma}$. 
Let $(A,[v])\in  \cU$ and
\begin{equation}\label{delav}
\delta_{\gamma}(A,[v])\colon T_{(A,[v])}\lagr
\times\PP(V)\lra \Sym^2(A\cap F_v)^{\vee}
\end{equation}
be as in~\eqref{difftre}. The restriction of $\delta_{\gamma}(A,[v])$ to the tangent space to $\lagr$ at $A$ is  surjective; thus  
\begin{equation}\label{libero}
\text{$\delta_{\gamma}(A,[v])$  is surjective.}
\end{equation}
Let $q\in\cX_{\gamma}$ and $f_{\cU}(q)=(A,[v])$. 
 Suppose that $q\notin\cW_{\cU}$ i.e.~that $\cork\gamma(p)\le 2$. If $\cork\gamma(p)=1$ then $Y_{\cU}=Y_{\gamma}$ is smooth because  the differential $\delta_{\gamma}(A,[v])$ is surjective: by~\Ref{prp}{gradodue} we get that $\cX_{\cU}$ is smooth at $q$.  If $\cork\gamma(p)=2$ then $\cX_{\cU}$ is smooth at $q$ by~\Ref{prp}{biancabernie} - recall that  the differential $\delta_{\gamma}(A,[v])$ is surjective. 
 This proves that $\sing\cX_{\cU}\subset\cW_{\cU}$. On the other hand $\cW_{\cU}\subset \sing\cX_{\cU}$ by~\Ref{clm}{sempresing}.
\end{proof}
We will close the present subsection by proving a few results about the individual $X_A$'s. 
\begin{lmm}\label{lmm:critliscio}
Let $A\in(\lagr\setminus\NN(V))$ and $[v]\in Y_A$. Suppose that 
$\dim(A\cap F_v)\le 2$ and that 
there is no $W\in\Theta_A$ (see~\eqref{eccoteta}) containing $v$.
Then $X_A$ is smooth at $f_A^{-1}([v])$.
\end{lmm}
\begin{proof}
Let $q\in f^{-1}_A([v])$.  Suppose that $\dim(A\cap F_v)=1$. By Corollary~2.5 of~\cite{og5} we get that $Y_A$ is smooth at $[v]$, thus $X_A$ is smooth at $q$  by~\Ref{prp}{gradodue}. 
Suppose that $\dim(A\cap F_v)=2$. Locally around $q$ the double cover $X_A\to Y_A$ is isomorphic to $X_{\ov{\gamma}}\to  Y_{\ov{\gamma}}$ where $\ov{\gamma}$ is  the symmetrization of the restriction of  $\beta_A$ to an affine neighoborhood $\Spec R$ of $[v]$. Thus we may consider the differential $\delta_{\ov{\gamma}}([v])$ - see~\eqref{difftre}. The differential is surjective by   Proposition~2.9 of~\cite{og5}, thus 
 $X_A$ is smooth at $q$ by~\Ref{prp}{biancabernie}. 
\end{proof}
\begin{prp}\label{prp:buono}
Let $A\in(\lagr\setminus\NN(V))$. Then  $X_A$ is smooth if and only if $A\in\lagr^0$. 
\end{prp}
\begin{proof}
If $A\in\lagr^0$ then $X_A$ is smooth by~\cite{og2}. Suppose that $X_A$ is smooth. Then $A\notin\Delta$ by~\Ref{clm}{sempresing}. Assume that $A\in\Sigma$; we will reach a contradiction. Let $W\in\Theta_A$ and $[v]\in\PP(W)$ - notice that $\PP(W)\subset Y_A$. Let $q\in f_A^{-1}([v])$. Since $A\notin\Delta$ we have $1\le \dim(A\cap F_v)\le 2$.  Suppose that $\dim(A\cap F_v)=1$. Then $Y_A$ is singular at $[v]$ by  Corollary~2.5 of~\cite{og5}, thus $X_A$ is singular at $q$ by~\Ref{prp}{gradodue}. Suppose that $\dim(A\cap F_v)=2$. Let $\ov{\gamma}$ be as in the proof of~\Ref{lmm}{critliscio}. Then $\delta_{\ov{\gamma}}([v])$ is not surjective - see Proposition 2.3 of~\cite{og5} - and hence 
 $X_A$ is singular at $q$ by~\Ref{prp}{biancabernie}.
\end{proof}
\subsection{The desingularization}\label{subsec:oradesing}
\setcounter{equation}{0}
\begin{dfn}\label{dfn:piazzata}
Let $\lagr^{*}\subset \lagr$ be the set of $A$ such that the following hold:
\begin{enumerate}
\item[(1)]
$A\notin\NN(V)$.
\item[(2)]
 $Y_A[3]$ is finite.
\item[(3)]
 $Y_A[3]=Y_A(3)$.
\end{enumerate}
\end{dfn}
\begin{rmk}\label{rmk:piazzata}
$\lagr^{*}$ is an open subset of $\lagr$.
\end{rmk}
\begin{clm}\label{clm:nodeco}
$(\lagr\setminus\Sigma)\subset\lagr^{*}$. 
\end{clm}
\begin{proof}
 Item~(1) of~\Ref{dfn}{piazzata} holds by Claim 2.11 of~\cite{og5}. Let's prove that Item~(2) of~\Ref{dfn}{piazzata} holds. Suppose that $Y_A[3]\not= Y_A(3)$ i.e.~there exists  $[v_0]\in\PP(V)$ such that $\dim(A\cap F_{v_0})\ge 4$. Let $V_0\subset V$ be a codimension-$1$ subspace transversal to $[v_0]$ and let $\rho^{v_0}_{V_0}$ be as in~\eqref{rocco}. Let 
${\bf K}:=\PP(\rho^{v_0}_{V_0}(A\cap F_{v_0}))$. Then $\dim{\bf K}\ge 3$; since $\Gr(2,V_0)$ has codimension $3$ in $\PP(\bigwedge^2 V_0)$ it follows that there exists $[\alpha]\in {\bf K}\cap\Gr(2,V_0)$. Let $\wt{\alpha}\in (A\cap F_{v_0})$ such that $\rho^{v_0}_{V_0}(\wt{\alpha})=\alpha$. Then $\wt{\alpha}$ is non-zero and decomposable, that is a contradiction because $A\notin\Sigma$. Lastly let's prove that Item~(3) of~\Ref{dfn}{piazzata} holds. Let $[v_0]\in Y_A[3]=Y_A(3)$. Then $(A,[v_0])\in \wt{\Delta}(0)$. Let $K:=A\cap F_{v_0}$ and $\tau^{v_0}_K$ be as in~\eqref{botte}. We have
\begin{equation*}
T_{[v_0]}Y_A[3]=T_{[v_0]}Y_A(3)=\ker\tau^{v_0}_K.
\end{equation*}
By~\Ref{prp}{aloha} the map  $\tau^{v_0}_K$ is injective. Thus $[v_0]$ is an isolated point of $Y_A[3]$. 
\end{proof}
Let $A\in\lagr^{*}$.  Let $\cU\subset\lagr^{*}$ be a small open (either in the Zariski or in the \emph{classical} topology) subset containing $A$.
In particular $\rho_{\cU}\colon \cX_{\cU}\to\cY_{\cU}$ exists. 
Let $\pi_{\cU}\colon\wt{\cX}_{\cU}\to\cX_{\cU}$ be the blow-up of $\cW_{\cU}$ and $E_{\cU}$ be the exceptional set of $\pi_{\cU}$. 
\begin{clm}\label{clm:scoppioliscio}
Keep notation as above. Then  $\wt{\cX}_{\cU}$ is smooth. If $\cU$ is open and sufficiently small in the \emph{classical} topology then  
we have a locally-trivial fibration    
\begin{equation}\label{tilli}
E_{\cU} \lra    Y_{\cU}[3].
\end{equation}
Let $(A,[v])\in Y_{\cU}[3]$. The fiber of~\eqref{tilli}  over $(A,[v])$ is  isomorphic to $\PP(A\cap F_v)^{\vee}\times\PP(A\cap F_v)^{\vee}$ and the restriction of  $N_{{E_{\cU}}/\wt{\cX}_{\cU}}$ to the fiber is isomorphic to $\cO_{\PP(A\cap F_v)^{\vee}}(-1)\boxtimes\cO_{\PP(A\cap F_v)^{\vee}}(-1)$.  
\end{clm}
\begin{proof}
By~\Ref{prp}{singfam} we know that $\wt{\cX}_{\cU}$ is smooth outside $E_\cU$. It remains to examine  $\wt{\cX}_{\cU}$ over $\cW_{\cU}\cong\cY_{\cU}[3]$. We may assume that $\cU=U_B\times\cN$ is as in the proof of~\Ref{prp}{singfam}. We will adopt the notation of that proof. Let $q\in\cX_{\gamma}$ and $f_{\cU}(q)=(A,[v])=p$. A neighborhood of $q$ in $X_{\cU}$ is isomorphic to $X_{\gamma}$ where $\gamma$ is given by~\eqref{gamfam} - see~\eqref{pavia}.  
 We are assuming that $q\in\cW_{\cU}$ and hence $\cork\gamma(p)=3$. Let $f\colon X(\cV)\to Y(\cV)$ be as in~\Ref{subsec}{doppipiccoli} i.e.~$f$ is the universal double covering of corank $3$ at the origin. We claim that there exists a map $\nu\colon X_{\gamma}\to X(\cV)$ such that the following diagram commutes
\begin{equation}\label{quadruniv}
\xymatrix{  X_{\gamma} \ar_{f_{\gamma}}[d]  \ar^{\nu}[r]  & 
X(\cV) \ar^{f}[d] \\   
 Y_{\gamma}  \ar^{\mu}[r] &  Y(\cV)}
\end{equation}
and  $X_{\gamma}$ is identified with the fibered product $Y_{\gamma}\times_{Y(\cV)}X(\cV)$.
  In fact it suffices to apply the reduction procedure of~\Ref{subsec}{prodenne} that leads to~\Ref{clm}{nucleo}.
 Let ${\bf K}$ be as in~\Ref{clm}{nucleo}: by~\eqref{yuguali} we have $(Y_{\gamma_{\bf K}},p)=
(Y_{\gamma},p)$ and by~\Ref{clm}{nucleo} we have a natural isomorphism $(X_{\gamma_{\bf K}},f_{\gamma_{\bf K}}^{-1}(p))\overset{\sim}{\to}
(X_{\gamma},f_{\gamma}^{-1}(p))$ commuting with $f_{\gamma_{\bf K}}$ and  $f_{\gamma}$. Let  $\cU=\Spec R$:  
we are free to replace $\cU$ by any affine open subset containing $(A,[v])$. 
Thus we may assume that ${\bf K}$ is a trivial $R$-module i.e.~${\bf K}=\cV\otimes R$ where $\cV$ is a complex $3$-dimensional vector-space. Hence we may view $\gamma_{\bf K}$ as a map $\gamma_{\bf K}\colon \Spec R\to \Sym^2\cV^{\vee}$. Notice that we have equality of schemes $Y_{\gamma}=\gamma_{\bf K}^{-1}Y(\cV)$; thus   the restriction of $\gamma_{\bf K}$ to $Y_{\gamma} $ defines a map $\mu\colon Y_{\gamma} \to  Y(\cV)$. The  claim follows.   By surjectivity of 
$\delta_{\gamma}(A,[v])$ - see~\eqref{libero} - we get that the germ $(X_{\gamma},f_{\gamma}^{-1}(p))$ is the product of a smooth germ (of dimension $54$) and  the germ $(X(\cV),f^{-1}(0))$. Looking at the explicit description of $X(\cV)$ given by~\Ref{prp}{equarela} we get right away that $\wt{\cX}_{\cU}$ is smooth over $q$ and the remaining statements as well. We need to assume that $\cU$ is a small open subset in the classical topology in order to ensure that Map~\eqref{tilli} is a locally-trvial fibration.
\end{proof}
\begin{rmk}\label{rmk:induista}
Let $A\in\lagr^{*}$ and let   $Y_A[3]=\{[v_1],\ldots,[v_s]\}$. 
Let $\cU\subset\lagr^{*}$ be a small open (in the \emph{classical} topology) subset containing $A$. 
For each $1\le i\le s$ choose a projection 
\begin{equation}\label{epsixa}
E_{\cU}([v_i])\lra  \PP(A\cap F_v)^{\vee}. 
\end{equation}
There exists a unique   $\PP^2$-fibration
\begin{equation}\label{marchette}
\epsilon\colon E_{\cU}\lra \star
\end{equation}
where $\star$ is itself a fibration over $Y_{\cU}[3]$ with fiber $\PP(A\cap F_v)^{\vee}$ over $(A,[v])$.  We say that~\eqref{epsixa} is a {\it choice of $\PP^2$-fibration $\epsilon$ for $X_A$}.
\end{rmk}
Let $A\in\lagr^{*}$ and choose a $\PP^2$-fibration $\epsilon$ for $X_A$. Let $\cU\subset\lagr^{*}$ be a small open (in the \emph{classical} topology) subset containing $A$.  By~\Ref{clm}{scoppioliscio} the normal bundle of $E_{\cU}$ along the fibers
 of~\eqref{marchette} is $\cO_{\PP^2}(-1)$. Thus there exists a contraction $c_{\cU,\epsilon}\colon\wt{\cX}_{\cU}\to 
\cX^{\epsilon}_{\cU}$  in the category of complex manifolds fitting into a commutative diagram
\begin{equation}\label{contraggo}
\xymatrix{ \wt{\cX}_{\cU} \ar_{\pi_{\cU}^{\epsilon}}[dr]  \ar^{c_{\cU,\epsilon}}[rr]   &  &
\cX^{\epsilon}_{\cU} \ar^{g_{\cU}^{\epsilon}}[dl] \\   
& \cX_{\cU} &}
\end{equation}
Let $f_{\cU}^{\epsilon}=f_{\cU}\circ g_{\cU}^{\epsilon}\colon \cX^{\epsilon}_{\cU}\to \cY_{\cU}$ and
$\rho_{\cU}^{\epsilon}\colon \cX_{\cU}^{\epsilon}\to \cU$ be the map $f_{\cU}^{\epsilon}$ followed by  $\cY_{\cU}\to\cU$. Let 
\begin{equation*}
X^{\epsilon}_A:=  (\rho^{\epsilon}_{\cU})^{-1}(A),\quad
g_A^{\epsilon}:=  g_{\cU}^{\epsilon}|_{X_A^{\epsilon}},\quad
f_A^{\epsilon}:=  f_{\cU}^{\epsilon}|_{X_A^{\epsilon}},\quad
\cO_{X_A^{\epsilon}}(1):= (f_A^{\epsilon})^{*}\cO_{Y_A}(1),\quad
H_A^{\epsilon}\in   |\cO_{X_A^{\epsilon}}(1)|.
\end{equation*}
Our notation does not make any reference to $\cU$ because the isomorphism class of the polarized couple $(X^{\epsilon}_A,\cO_{X_A^{\epsilon}}(1))$ does not depend on the open set $\cU$ containing $A$. Notice that if $A\in\Delta $ then $\cO_{X_A^{\epsilon}}(1)$ is not ample, in fact it is trivial on  $s$ copies of $\PP^2$ where $s=|Y_A[3]|$.
Of course 
\begin{equation}
\text{$(X^{\epsilon}_A,\cO_{X_A^{\epsilon}}(1))\cong
(X_A,\cO_{X_A}(1))$ if $A\in(\lagr\setminus\Delta)$.} 
\end{equation}
\begin{prp}\label{prp:fuorisig}
Let $A\in\lagr^{*}$  and let ${\epsilon}$ be a choice of $\PP^2$-fibration for $X_A$. 
\begin{itemize}
\item[(1)]
$X_A^{\epsilon}$  is smooth away from $(f_A^{\epsilon})^{-1}(\bigcup_{W\in\Theta_A}\PP(W))$.
\item[(2)]
If $[v_i]\in Y_A[3]$ then
$(f_A^{\epsilon})^{-1}[v_i]\cong \PP(A\cap F_{v_i})^{\vee}$.
\item[(3)]
If $\epsilon'$ is another choice of $\PP^2$-fibration for $X_A$
 there exists a commutative diagram
\begin{equation}
\xymatrix{ X_A^{\epsilon} \ar_{f_A^{\epsilon}}[dr]   & \dra &
X_A^{\epsilon'} \ar^{f_A^{\epsilon'}}[dl] \\   
& Y_A &}
\end{equation}
where the birational map is the flop of a collection of $(f_A^{\epsilon})^{-1}[v_i]$'s. Conversely every flop of a collection of $(f_A^{\epsilon})^{-1}[v_i]$'s is isomorphic to one  $X_A^{\epsilon'}$.
\end{itemize}
\end{prp}
\begin{proof}
Let's prove Item~(1). $X_A^{\epsilon}$  is smooth away from $(f_A^{\epsilon})^{-1}(Y_A[3]\cup\bigcup_{W\in\Theta_A}\PP(W))$ by~\Ref{lmm}{critliscio}. It remains to prove that  $X_A^{\epsilon}$  is smooth at every point of  $(f_A^{\epsilon})^{-1}\{[v_1],\ldots,[v_s]\}$ where 
\begin{equation}
\{[v_1],\ldots,[v_s]\}= Y_A[3]\setminus\bigcup_{W\in\Theta_A}\PP(W).
\end{equation}
Let $\cU\subset\lagr^{*}$ be a small open (in the \emph{classical} topology) subset containing $A$. 
  Let $\wt{\rho}_{\cU}:=\rho_{\cU}\circ\pi_{\cU}$; thus $\wt{\rho}_{\cU}\colon\wt{X}_{\cU}\to\cU$. 
For $1\le i\le s$ the fiber over $(A,[v_i])$ of Fibration~\eqref{tilli} is canonically isomorphic to $\PP(A\cap F_{v_i})^{\vee}\times\PP(A\cap F_{v_i})^{\vee}$. Let $\wh{X}_A\subset\wt{X}_{\cU}$ be the strict transform of $X_A$. Abusing notation we write
\begin{equation}
\wt{\rho}_{\cU}^{-1}(A)=\wh{X}_A\cup
\bigcup\limits_{i=1}^s \PP(A\cap F_{v_i})^{\vee}\times\PP(A\cap F_{v_i})^{\vee}\,.
\end{equation}
(Of course $\PP(A\cap F_{v_i})^{\vee}\times\PP(A\cap F_{v_i})^{\vee}$ denotes the fiber over $(A,[v_i])$ of Fibration~\eqref{tilli}.)
The components  $\PP(A\cap F_{v_i})^{\vee}\times\PP(A\cap F_{v_i})^{\vee}$ are pairwise disjoint. We claim that for $i=1,\ldots,s$ the intersection 
\begin{equation}
E_{A,i}:=\wh{X}_A\cap(\PP(A\cap F_{v_i})^{\vee}\times\PP(A\cap F_{v_i})^{\vee})
\end{equation}
is a smooth symmetric divisor in the linear system $|\cO_{\PP(A\cap F_{v_i})^{\vee}}(1)\boxtimes\cO_{\PP(A\cap F_{v_i})^{\vee}}(1)|$. In order to prove this we go back to Map~\eqref{mappasigma} - recall that $\cV$ is a $3$-dimensional complex vector space. Pull-back by $\sigma$ defines an isomorphism 
\begin{equation}
\Sym^2\cV^{\vee}\overset{\sigma^{*}}{\lra}
(\cV^{\vee}\otimes\cV^{\vee})^{\ZZ/(2)}=:\Sym_2\cV^{\vee}
\end{equation}
which is $\GL(\cV)$-equivariant. Isomorphism $\sigma^{*}$ induces a $\PGL(\cV)$-equivariant isomorphism of projective spaces ${\bf p}\colon \PP(\Sym^2\cV^{\vee})\overset{\sim}{\lra}
\PP(\Sym_2\cV^{\vee})$.
Of course ${\bf p}$ maps a point in the unique open $\PGL(\cV)$-orbit  
of $\PP(\Sym^2\cV^{\vee})$ to a point in the unique open $\PGL(\cV)$-orbit  
of $\PP(\Sym_2\cV^{\vee})$. Now let $\cV=(A\cap F_{v_i})^{\vee}$. Let $K_i:=(A\cap F_{v_i})$ and $\tau^{v_i}_{K_i}$ be as in \eqref{botte}.
By~\Ref{prp}{aloha}  we have that $\im(\tau^{v_i}_{K_i})$   belongs to the unique open $\PGL(K_i)$-orbit  
of $\PP(\Sym^2(A\cap F_{v_i}))$. Commutative Diagram~\eqref{triangolo} gives that $E_{A,i}$ is a symmetric smooth divisor in $|\cO_{\PP(A\cap F_{v_i})^{\vee}}(1)\boxtimes\cO_{\PP(A\cap F_{v_i})^{\vee}}(1)|$. Thus we have described $\wt{\rho}_{\cU}^{-1}(A)$. Since $X^{\epsilon}_{\cU}$ is obtained from $\wt{X}_{\cU}$ by  contracting  $E_{\cU}$ along the $\PP^2$-fibration $\epsilon$ it follows that $X_A^{\epsilon}$  is smooth at every point of  $(f_A^{\epsilon})^{-1}\{[v_1],\ldots,[v_s]\}$.
 This proves Item~(1). Since $X^{\epsilon}_A$ is obtained from $\wh{X}_A$ by contracting each of the divisors $E_{A,i}$ along the fibration $\PP^1\to E_{A,i}\to  \PP(A\cap F_{v_i})^{\vee}$ determined by $\epsilon$ (and similarly for $\epsilon'$) we also get Items~(2) and~(3). 
\end{proof}
\begin{crl}
Let $A\in(\lagr\setminus\Sigma)$. Then   $g_A^{\epsilon}\colon X_A^{\epsilon}\to X_A$ is a desingularization for every choice of $\PP^2$-fibration $\epsilon$ for $X_A$. 
\end{crl}
\begin{proof}
By~\Ref{clm}{nodeco} we know that $A\in\lagr^{*}$: thus~\Ref{prp}{fuorisig} applies to $X_A^{\epsilon}$. Since $A\notin\Sigma$ we get that  $X_A^{\epsilon}$ is smooth by Item~(1) of~\Ref{prp}{fuorisig}. 
\end{proof}
\begin{crl}\label{crl:unadef}
Let $A,A'\in(\lagr\setminus\Sigma)$ and  ${\epsilon},\epsilon'$ be choices of $\PP^2$-fibration for $X_A$. The quasi-polarized $4$-folds $(X_{A}^{\epsilon},H_{A}^{\epsilon})$ and $(X_{A'}^{\epsilon},H_{A'}^{\epsilon})$ are deformation equivalent.
\end{crl}
 \section{Double EPW-sextics parametrized by $\Delta$}\label{sec:zitelle}
 \setcounter{equation}{0}  
Let $A\in\Delta$ and $[v_0]\in Y_A(3)$. In the first subsection we will associate to $(A,[v_0])$ (under some  hypotheses which are certainly satisfied if $A\notin\Sigma$) a $K3$ surface $S_A(v_0)$ of genus $6$,  meaning that it comes equipped with a big and nef divisor class $D_A(v_0)$ of square $10$. We will also prove a converse:  given a generic such pseudo-polarized $K3$ surface $S$ there exist $A\in\Delta$ and $[v_0]\in Y_A(3)$ such that the pseudo-polarized surfaces $S$ and  $S_A(v_0)$ are isomorphic.  In the second subsection we will assume that $A\in(\Delta\setminus\Sigma)$ - with this hypothesis   $D_A(v_0)$ is very ample. We will prove that there exists a bimeromorphic map $\psi\colon S_A^{[2]}(v_0)\dashrightarrow X_A^{\epsilon}$ where $\epsilon$ is an arbitrary choice of $\PP^2$-fibration for $X_A$. That such a map exists for generic $A\in\Delta$ could be proved by invoking the results of~\cite{og4}.  Here we will present a direct proof (we will not appeal to~\cite{og4} nor to~\cite{og2}).
 Moreover we will prove that if  $S_A(v_0)$ contains no lines  (this will be the case for generic $A$)  then there exists a choice of $\epsilon$ for which
 $\psi$ is regular - in particular $X_A^{\epsilon}$ is projective for such $\epsilon$. Lastly we will notice that the above results show that a smooth double cover of an EPW-sextic is a deformation of the Hilbert square of a $K3$ (and that the family of double EPW-sextics is a locally versal family of projective Hyperk\"ahler manifolds): the proof is more direct than the proof of~\cite{og2}. 
\subsection{EPW-sextics and $K3$ surfaces}\label{subsec:kapdel}
\setcounter{equation}{0}
\begin{ass}\label{ass:avuzero}
$A\in\lagr$, $[v_0]\in Y_A(3)$ and the following hold:
\begin{itemize}
\item[(a)]
There exists a codimension-$1$ subspace $V_0\subset V$ such that $\bigwedge^3 V_0\pitchfork A$ i.e.~$\bigwedge^3 V_0\cap A=\{0\}$. 
\item[(b)]
There exists at most 
 one $W\in\Theta_A$ containing $v_0$.
\item[(c)]
If $W\in\Theta_A$ contains $v_0$ then $A\cap (\bigwedge^2 W\wedge V)=\bigwedge^3 W$.
\end{itemize}
\end{ass}
\begin{rmk}\label{rmk:pixmania}
Let $A\in(\Delta\setminus\Sigma)$. Let $[v_0]\in Y_A(3)$ ($=Y_A[3]$ by~\Ref{clm}{nodeco}). Then~\Ref{ass}{avuzero} holds. In fact Items~(b) and~(c) hold trivially while Item~(a) holds by Claim~2.11 and Equation~(2.81) of~\cite{og5}. 
\end{rmk}
Let $(A,[v_0])$ be as in~\Ref{ass}{avuzero}: we will  define a surface $S_A(v_0)$ of genus $6$. 
The condition that $\bigwedge^3 V_0$ is transverse to $A$ is open: thus we may assume that we have a direct-sum decomposition
\begin{equation}\label{trasuno}
V=[v_0]\oplus V_0.
\end{equation}
We will denote by $\cD$ be the direct-sum decomposition of $V$ appearing in~\eqref{trasuno}. 
Let 
\begin{equation}\label{kappagrasso}
K^{\cD}_A:= \rho^{v_0}_{V_0}(A\cap F_{v_0}).
\end{equation}
where $\rho^{v_0}_{V_0}$ is given by~\eqref{rocco}. 
Choose a volume-form on $V_0$. Wedge-product followed by the volume-form defines an isomorphism $\bigwedge^3 V_0\cong\bigwedge^2 V_0^{\vee}$ and hence it makes sense to let 
\begin{equation}
F^{\cD}_A:=\PP(\Ann  K^{\cD}_A)\cap\Gr(3,V_0).
\end{equation}
By~\Ref{prp}{fanoindue} and~\Ref{prp}{singfano} (see the Appendix) we know that $F^{\cD}_A$ is  a Fano $3$-fold with at most one singular point. Next we will define a  quadratic form  on $\Ann  K^{\cD}_A$.  By Item~(a) of~\Ref{ass}{avuzero} the subspace $A$ is the graph of a map $\wt{q}_A^{\cD}\colon \bigwedge^2 V_0\to  \bigwedge^3 V_0$: explicitly 
\begin{equation}\label{giannibrera}
\wt{q}_A^{\cD}(\alpha)=\beta\iff
(v_0\wedge\alpha+\beta)\in A. 
\end{equation}
The map $\wt{q}_A^{\cD}$  is symmetric because $A$, $\bigwedge^2 V_0$ and $\bigwedge^3 V_0$ are lagrangian subspaces of $\bigwedge^3 V$. Clearly $\ker\wt{q}_A^{\cD}=K^{\cD}_A$: it follows that  
 $\wt{q}_A^{\cD}$ induces an isomorphism
\begin{equation}
\wt{r}_A^{\cD}\colon \bigwedge^2 V_0/K^{\cD}_A\overset{\sim}{\lra} 
\Ann  K^{\cD}_A\subset \bigwedge^3 V_0.
\end{equation}
The inverse $(\wt{r}_A^{\cD})^{-1}$ defines a non-degenerate quadratic form $(r_A^{\cD})^{\vee}$ on $\Ann  K^{\cD}_A$. For future reference we unwind the definition of $(\wt{r}_A^{\cD})^{-1}$ and $(r_A^{\cD})^{\vee}$. 
 Let $\beta\in \Ann   K^{\cD}_A$ i.e.
\begin{equation}\label{killk}
v_0\wedge\alpha+\beta\in A,\qquad \alpha\in\bigwedge^2 V_0.
\end{equation}
 Then
\begin{equation}\label{errepol}
(\wt{r}_A^{\cD})^{-1}(\beta)\equiv  \alpha\pmod{K^{\cD}_A},\qquad 
(r_A^{\cD})^{\vee}(\beta)=  \vol(v_0\wedge\alpha\wedge\beta).
\end{equation}
Let $V((r_A^{\cD})^{\vee})\subset\PP(\Ann   K^{\cD}_A)$ be the zero-scheme of $(r_A^{\cD})^{\vee}$:  a smooth $5$-dimensional quadric. Let
 \begin{equation}\label{essdec}
S_A^{\cD}:=V((r_A^{\cD})^{\vee})\cap F^{\cD}_A.
\end{equation}
Our first goal is to show that  $S_A^{\cD}$ does not depend on the choice of the subspace $V_0\subset V$ complementary to $[v_0]$ i.e.~it depends only on $A$ and $[v_0]$. First we notice that $F^{\cD}_A$ is independent of $V_0$. In fact $\bigwedge^3 V_0$ is transversal to $F_{v_0}$; since both $\bigwedge^3 V_0$ and $F_{v_0}$ are Lagrangians the volume $\vol$ induces an isomorphism 
\begin{equation}\label{elipattino}
g_{V_0}\colon \bigwedge^3 V_0\overset{\sim}{\lra} F_{v_0}^{\vee}\,.
\end{equation}
Thus $g_{V_0}$ defines an inclusion
\begin{equation}\label{osiris}
 F^{\cD}_A\hra\PP(\Ann   K_A)\,.
\end{equation}
\begin{rmk}\label{rmk:zoppas}
The image of Map~\eqref{osiris} does not depend on $V_0$ i.e.~it depends exclusively on $A$ and $[v_0]\in Y_A(3)$; we  will denote it by $Z_A(v_0)$. 
\end{rmk}
Similarly $g_{V_0}$ defines an inclusion
\begin{equation}\label{dovesta}
{\bf g}_{V_0}\colon S_A^{\cD}\hra\PP(\Ann   K_A)\,.
\end{equation}
\begin{lmm}\label{lmm:intrinseco}
Keep notation and assumptions as above. Then ${\bf g}_{V_0}(S_A^{\cD})$  is independent of $V_0$, in other words it depends exclusively on $A$ and $[v_0]\in Y_A(3)$.
\end{lmm}
\begin{proof}
Let $V_0'\subset V$ be a codimension-$1$ subspace complementary to $[v_0]$ and transverse to $A$. 
Let $\cD'$ denote the corresponding direct-sum decomposition of $V$; we must show that  
\begin{equation}\label{stessosub}
{\bf g}_{V_0}(S_A^{\cD})={\bf g}_{V'_0}(S_A^{\cD'})\,.
\end{equation}
The subspace $V_0'$ is the graph of a linear function 
\begin{equation}
\begin{matrix}
V_0 & \lra & [v_0] \\
v & \mapsto & f(v)v_0
\end{matrix}
\end{equation}
 and hence we have an isomorphism
\begin{equation}
\begin{matrix}
V_0 & \overset{\psi}{\lra} & V_0' \\
v & \mapsto & v+f(v)v_0\,.
\end{matrix}
\end{equation}
We notice that
\begin{equation}\label{contrazione}
\bigwedge^ 3\psi(\beta)=\beta+v_0\wedge (f\contr\beta)
\end{equation}
where $\contr$ denotes contraction.
In particular $g_{V'_0}\circ\bigwedge^ 3\psi=g_{V_0}$. Moreover $\phi:=\bigwedge^ 3\psi|_{\Ann   K^{\cD}_A}$ is an isomorphism between $\Ann   K^{\cD}_A\subset\bigwedge^3 V_0$ and $\Ann   K^{\cD'}_{A'}\subset\bigwedge^3 V'_0$. 
Thus it suffices  to prove that 
\begin{equation}\label{corrisponde}
\phi(S_A^{\cD})=S_A^{\cD'}\,.
\end{equation}
We claim that
\begin{equation}\label{amenodi}
\phi^{*}(r_A^{\cD'})^{\vee}-(r_A^{\cD})^{\vee}
\in H^0(\cI_{F^{\cD}_A}(2))\,.
\end{equation}
In fact let $\beta\in \Ann   K^{\cD}_A\subset\bigwedge^3 V_0$; then~\eqref{killk} holds.
By~\eqref{contrazione} we get that
\begin{equation}
v_0\wedge(\alpha-(f\contr \beta))+\phi(\beta)=
v_0\wedge\alpha+\beta\in A\,.
\end{equation}
By~\eqref{contrazione} we get that
\begin{multline}
\phi^{*}(r_A^{\cD'})^{\vee}(\beta)=
\vol(v_0\wedge(\alpha-(f\contr\beta))\wedge\phi(\beta))=\\
\vol(v_0\wedge\alpha\wedge\phi(\beta))-
\vol(v_0\wedge (f\contr\beta)\wedge\phi(\beta))=\\
\vol(v_0\wedge\alpha\wedge\beta)-
\vol(v_0\wedge (f\contr\beta)\wedge\beta)=\\
(r_A^{\cD})^{\vee}(\beta)-
\vol(v_0\wedge (f\contr\beta)\wedge\beta)\,.
\end{multline}
The second term in the last expression  is  the restriction to $\PP(\Ann   K^{\cD}_A)$ of a Pl\"ucker quadratic form and hence it vanishes on $F^{\cD}_A$. This proves~(\ref{amenodi}) and hence~(\ref{corrisponde}) holds.  
\end{proof}
By the above lemma we may give the following definition.
\begin{dfn}
Let $A\in\lagr$. Suppose that $[v_0]\in Y_A(3)$ and that~\Ref{ass}{avuzero} holds. Let  $\cD$ be the  direct-sum decomposition~\eqref{trasuno}. We set
\begin{equation}
S_A(v_0):={\bf g}_{V_0}(S_A^{\cD}).
\end{equation}
\end{dfn}
Keep assumptions and notation as above. We single out special points of $S_A(v_0)$ as follows. Suppose that $W\in\Theta_A$ (see~\eqref{eccoteta} for the definition of $\Theta_A$) and assume that $v_0\notin W$. Let $\gamma$ be a generator of $\bigwedge^3 W$ i.e.~$\gamma$ is decomposable with $\supp(\gamma)=W$. By hypothesis $\bigwedge^3 V_0\cap A=\{0\}$ and hence $W\not\subset V_0$; thus 
\begin{equation}\label{mafalda}
\gamma=(v_0+u_1)\wedge u_2\wedge u_3,\qquad u_i\in V_0\,.
\end{equation}
Since $v_0\notin W$ we have $u_1\wedge u_2\wedge u_3\not=0$; thus $[u_1\wedge u_2\wedge u_3]\in F^{\cD}_A$. Moreover $[u_1\wedge u_2\wedge u_3]\in V((r_A^{\cD})^{\vee})$ by~\eqref{errepol} and hence $[u_1\wedge u_2\wedge u_3]\in S_A^{\cD}$. We let
\begin{equation}
\begin{matrix}
\Theta_A\setminus\{W\mid v_0\in W\} & 
\overset{\theta_A^{\cD}}{\lra} & S_A^{\cD} \\
W & \mapsto & [u_1\wedge u_2\wedge u_3]\,.
\end{matrix}
\end{equation}
The map
\begin{equation}
\theta_A(v_0):={\bf g}_{V_0}\circ \theta_A^{\cD}\colon (\Theta_A\setminus\{W\mid v_0\in W\})\to S_A(v_0)
\end{equation}
  is independent of $\cD$, i.e.~it depends exclusively on $A$ and $[v_0]$. Notice that $\theta_A(v_0)$ is injective. 
 \begin{prp}\label{prp:transex}
 Let $A\in\lagr$. Suppose that $[v_0]\in Y_A(3)$ and that~\Ref{ass}{avuzero} holds. Let  $\cD$ be the  direct-sum decomposition~\eqref{trasuno}. 
The set of points at which the intersection $V((r_A^{\cD})^{\vee})\cap F^{\cD}_A$ is not transverse is equal to
\begin{equation}
\im\theta_A^{\cD}\coprod(S_A^{\cD}\cap \sing F^{\cD}_A).
\end{equation}
\end{prp}
\begin{proof}
  Let $[\beta]\in S^{\cD}_A$. In particular $\beta$  is non-zero decomposable; let $U:=\supp \beta$. 
Moreover since $[\beta]\in F^{\cD}_A$ we have that~\eqref{killk} holds; let $\alpha\in\bigwedge^2 V_0$ be as in~\eqref{killk}.
We claim that 
\begin{equation}\label{vendette}
\text{$V((r^{\cD}_A)^{\vee})\pitchfork F^{\cD}_A$  at $[\beta]$  unless $\la\alpha,\,K^{\cD}_A\ra\cap\bigwedge^2 U\not=\es$.}
\end{equation}
In fact
the projective tangent space to $\Gr(3,V_0)$ at $[\beta]$ is given by
\begin{equation}\label{tangrass}
{\bf T}_{[\beta]}\Gr(3,V_0)=\PP(\Ann  (\bigwedge^2 U))\,.
\end{equation}
On the other hand~\eqref{errepol} gives that
\begin{equation}\label{tanaka}
{\bf T}_{[\beta]}V((r^{\cD}_A)^{\vee})=\PP(\Ann   \alpha)\cap\PP(\Ann   K^{\cD}_A)\,.
\end{equation}
Statement~(\ref{vendette}) follows at once from~(\ref{tangrass}) and~(\ref{tanaka}). 
Next we prove  that
\begin{equation}\label{sinequa}
\text{$\la\alpha,\,K^{\cD}_A\ra\cap\bigwedge^2 U\not=\es$ if and only if 
$[\beta]\in \sing F^{\cD}_A$ or $[\beta]\in \im\theta_A^{\cD}$.}
\end{equation}
Suppose that $[\beta]\in \sing F^{\cD}_A$; then Item~(1) of~\Ref{prp}{singfano} gives that $K^{\cD}_A\cap\bigwedge^2 U\not=\es$. 
Next suppose that $[\beta]\in \im\theta_A^{\cD}$; then $\alpha\in\bigwedge^2 U$ by~(\ref{mafalda}). 
This proves the \lq\lq if\rq\rq implication of~\eqref{sinequa}. Let us prove the \lq\lq only if\rq\rq implication. First assume that $K^{\cD}_A\cap\bigwedge^2 U\not=\{0\}$. Let $0\not=\kappa_0\in K^{\cD}_A\cap\bigwedge^2 U$. Then $\kappa_0$ is decomposable because $\dim U=3$ and hence $[\kappa_0]$ is the unique point belonging to $\PP(K^{\cD}_A)\cap \Gr(2,V_0)$. We get that $[\beta]$ is the unique singular point of $F^{\cD}_A$ by~\eqref{nelpiano}.  Lastly assume that $K^{\cD}_A\cap\bigwedge^2 U=\{0\}$. Then
there exists $\kappa\in K^{\cD}_A$ such that $(\alpha+\kappa)\in\bigwedge^2 U$. Since $\kappa\in K^{\cD}_A$ we have $(v_0\wedge(\alpha+\kappa)+\beta)\in A$. The tensor $(v_0\wedge (\alpha+\kappa)+\beta)\in A$ is  decomposable, let $W$ be its support. Then $v_0\notin W$ because  $\beta\not=0$ and hence $[\beta]=\theta_A^{\cD}(W)$. This finishes the proof of~\eqref{sinequa} and of the proposition.
\end{proof}
 \begin{crl}\label{crl:prekappa}
 Let $A\in\lagr$. Suppose that $[v_0]\in Y_A(3)$ and that~\Ref{ass}{avuzero} holds. Asssume in addition that   $\Theta_A$ is finite. 
Then $S_A(v_0)$ is a reduced and irreducible surface with
\begin{equation}\label{puntising}
\sing S_A(v_0)=
\im\theta_A(v_0)\coprod(S_A(v_0)\cap \sing Z_A(v_0))\,.
\end{equation}
(See~\Ref{rmk}{zoppas} for the definition of $Z_A(v_0)$.)
\end{crl}
\begin{proof}
By~\Ref{prp}{transex} we know that $S_A^{\cD}$ is a smooth surface outside the right-hand side of~\eqref{puntising}. By hypothesis   $\Theta_A$ is finite and hence the right-hand side of~\eqref{puntising} is finite. On the other hand by~\Ref{prp}{singfano} we know that $Z_A(v_0)$ is a $3$-fold with at most one singular point, necessarily an ordinary quadratic singularity, and   $S_A^{\cD}$ is the complete intersection of $Z_A(v_0)$ and a quadric hypersurface. 
 It follows that $S_A^{\cD}$ is reduced and irreducible with singular set as claimed.
\end{proof}
 \begin{crl}\label{crl:kappatre}
Let hypotheses be as in~\Ref{crl}{prekappa}. Suppose in addition that  $S_A(v_0)$ has Du Val singularities. Let $\wh{S}_A(v_0)\to S_A(v_0)$ be the minimal desingularization. Then $\wh{S}_A(v_0)$ is  a  $K3$ surface.
\end{crl}
\begin{proof}
Let $\cO_{Z_A(v_0)}(1)$ be the pull-back by Map~(\ref{osiris}) of the hyperplane line-bundle  on $\PP(\Ann   (F_{v_0}\cap A))$. Then $S_A(v_0)\in |\cO_{Z_A(v_0)}(2)|$. By~\Ref{prp}{fanoindue} and~\Ref{prp}{singfano} there exist smooth divisors in $|\cO_{Z_A(v_0)}(2)|$ and they are $K3$ surfaces; by simultaneous resolution of Du Val singularities we get that $\wh{S}_A(v_0)$ is  a  $K3$ surface.
\end{proof}
\begin{crl}\label{crl:buonogrande}
Let $A\in(\Delta\setminus\Sigma)$. Let $[v_0]\in Y_A(3)$ (and hence~\Ref{ass}{avuzero} holds by~\Ref{rmk}{pixmania}). Then $S_A(v_0)$ is a (smooth) $K3$.
\end{crl}
\begin{proof}
Immediate consequence of~\Ref{crl}{kappatre}.  
\end{proof}
Under the hypotheses of~\Ref{crl}{kappatre} let $\cO_{S_A(v_0)}(1)$ be the restriction to $S_A(v_0)$ of $\cO_{Z_A(v_0)}(1)$. Let   $\cO_{\wh{S}_A(v_0)}(1)$ be the pull-back of  $\cO_{S_A(v_0)}(1)$ to $\wh{S}_A(v_0)$. We set
\begin{equation}
D_A(v_0)\in  | \cO_{S_A(v_0)}(1) | \qquad
\wh{D}_A(v_0)\in  | \cO_{\wh{S}_A(v_0)}(1) |.
\end{equation}
\begin{rmk}\label{rmk:genersei}
Let hypotheses be as in~\Ref{crl}{kappatre}. Then $(\wh{S}_A(v_0),\wh{D}_A(v_0))$ is a quasi-polarized 
$K3$ surface of genus $6$. Moreover the composition
\begin{equation}
\wh{S}_A(v_0) \lra S_A(v_0)\lra \PP(\Ann   (F_{v_0}\cap A))
\end{equation}
is identified (up to projectivities) with the map associated to the complete linear system $| \wh{D}_A(v_0) |$.
\end{rmk}
\Ref{rmk}{genersei} has a converse; in order to formulate it 
we identify $F_{v_0}\cong \bigwedge^2 (V/[v_0])$ (the identification is well-defined up to homothety).
 \begin{ass}\label{ass:kappass}
 $ K\in\Gr(3, F_{v_0})$ and
\begin{enumerate}
\item[(1)]
$\PP(K)\cap \Gr(2,V/[v_0])=\es$, or
\item[(2)]
the scheme-theoretic intersection $\PP(K)\cap \Gr(2,V/[v_0])$ is a single reduced point.
\end{enumerate}
 \end{ass}
 Let 
\begin{equation}\label{voodoo}
W_K:=\PP(\Ann   K)\cap \Gr(3,V/[v_0]).
\end{equation}
 (This makes sense because we have an isomorphism $\bigwedge^2 (V/[v_0])\overset{\sim}{\lra}\bigwedge^3 (V/[v_0])^{\vee}$ well-defined up to homothety). Let
\begin{equation}\label{bombay}
S:=W_K\cap Q,\qquad \text{$Q\subset\PP(\Ann   K)$ a quadric.}
\end{equation}
If $Q$ is generic then $S$ is a linearly normal $K3$ surface of genus $6$, see~\Ref{crl}{kappatre}. In fact the family of such $K3$ surfaces is locally versal. More generally suppose that~\Ref{ass}{kappass} holds, that $S$ is given by~\eqref{bombay} and that  $S$
has DuVal singularities. Let  $\wh{S}\to S$ be the minimal desingularization - thus $\wh{S}$ is a $K3$ surface.   Let $D\in|\cO_S(1)|$ and $\wh{D}$ be the pull-back of $D$ to $\wh{S}$. 
Consider the family $\cS\to B$ of deformations  of $(S,D)$  obtained by deforming slightly $K$ and $Q$;
by Brieskorn and Tjurina there is a suitable base change $\wh{B}\to B$ such that the pull-back of $\cS$ to $\wh{B}$ admits a simultaneous resolution of singularities $\wh{S}\to\wh{B}$ with fiber $\wh{S}$ over the point corresponding to $S$. Of course there is a divisor class $\wh{\cD}$ on $\wh{\cS}$ whose restriction to $\wh{S}$ is $\wh{D}$ - thus $\wh{\cS}\to\wh{B}$ is a family of quasi-polarized $K3$ surfaces.  The following result is well-known  - we omit the (standard) proof.
\begin{prp}\label{prp:generikappa}
Keep notation and hypotheses as above. 
The family $\wh{\cS}\to\wh{B}$ is a versal family of quasi-polarized $K3$ surfaces.  
\end{prp}
\begin{lmm}\label{lmm:scabrosa}
Suppose that~\Ref{ass}{kappass} holds. Let $S$ be as in~\eqref{bombay} and assume that $Q$ is transversal to $W_K$ outside a finite set - thus $S$ is a surface  with  finite singular set.  There exists a smooth quadric $Q'\subset\PP(\Ann   K)$  such that $S=W_K\cap Q'$.
\end{lmm}
\begin{proof}
Since $W_K$ is cut out by quadrics Bertini's Theorem gives that the generic quadric in $\PP(\Ann   K)$ containing $S$ is smooth outside $\sing S$; let $Q_0=V(P_0)$ be such a  quadric. Let $p\in \sing S$. The generic quadric $Q'=V(P')\in | \cI_{W_K}(2) |$ is smooth at $p$ and hence $V(P_0+P')$ is smooth at $p$. Since $\sing S$ is finite we get that the generic quadric $Q$ containing $S$  is smooth at all points of $\sing S$.  It follows that the generic quadric $Q$ containing $S$ is smooth. 
\end{proof}
The following corollary provides an inverse of the process which produces $S_A(v_0)$ out of  $(A,[v_0])\in\wt{\Delta}(0)$ (with the extra hypotheses in~\Ref{ass}{avuzero}). 
\begin{prp}\label{prp:pellegrini}
Suppose that~\Ref{ass}{kappass} holds.  Let $S$ be as in~\eqref{bombay} and assume that
$Q$ is smooth and transversal to $W_K$ outside a finite set.
There exist $A\in\Delta$, $[v_0]\in\PP(V)$ and a codimension-$1$ subspace $V_0\subset V$ transversal to $[v_0]$ such that the following hold:
\begin{itemize}
\item[(1)]
$\bigwedge^3 V_0\cap A=\{0\}$,
\item[(2)]
Items~(c) and~(d) of~\Ref{ass}{avuzero} hold,
\item[(3)]
the natural isomorphism $\PP(\bigwedge^3(V/[v_0]))\overset{\sim}{\lra}\PP(\bigwedge^3 V_0)$ maps $S$ to $S_A^{\cD}$ where $\cD$ is the direct-sum decomposition of $V$ appearing in~\eqref{trasuno}.  
\end{itemize}
If we replace the quadric $Q$ by a smooth quadric $Q'\subset\PP(\Ann   K)$ such that $S=W_K\cap Q'$ and let $A'\in\Delta$ be the corresponding point, there exists a projectivity of $\PP(V)$ fixing $[v_0]$ which takes $A$ to $A'$.  
\end{prp}
\begin{proof}
Let $Q=V(P)$. The dual of $\Ann   K$ is $\bigwedge^2(V/[v_0])/K$; thus the polarization of $P$ defines a non-degenerate symmetric  
map
\begin{equation}
\Ann   K\overset{\sim}{\lra} \bigwedge^2(V/[v_0])/K.
\end{equation}
The inverse of the above map is 
non-degenerate symmetric  
map
\begin{equation}
\bigwedge^2(V/[v_0])/K\overset{\sim}{\lra}\Ann   K. 
\end{equation}
Composing on the right with $\bigwedge^2(V/[v_0])\overset{\sim}{\lra}\bigwedge^2(V/[v_0])$ and the quotient map $\bigwedge^2(V/[v_0])\to \bigwedge^2(V/[v_0])/K$ and on the left with $\Ann   K\hra\bigwedge^3(V/[v_0])$ and  $\bigwedge^3(V/[v_0])\overset{\sim}{\lra}\bigwedge^3(V/[v_0])$ we get a symmetric map
\begin{equation}
\bigwedge^2 V_0\lra \bigwedge^3 V_0
\end{equation}
with $3$-dimensional kernel corresponding to $K$. The graph of the above map is a Lagrangian $A\in\lagr$. One checks easily that~(1), (2) and~(3) hold. One gets that the projective equivalence of $A$ does not depend on $Q$ by going through the proof of~\Ref{lmm}{intrinseco}.  
\end{proof}
\subsection{$X_A^{\epsilon}$ for $A\in(\Delta\setminus\Sigma)$}\label{subsec:hilbdue}
\setcounter{equation}{0}

 Let  $S$ be a $K3$. 
 Let $\Delta_S^{[2]}\subset S^{[2]}$ be the irreducible codimension $1$ subset parametrizing non-reduced subschemes. There exists  a square root of the line bundle $\cO_{S^{[2]}}(\Delta_S^{[2]})$: we denote by $\xi$ its first Chern class. 
 There is a natural  morphism of integral Hodge structures $\mu\colon H^2(S)\to H^2(S^{[2]})$ such that 
 $H^2(S^{[2]};\ZZ)=\mu(H^2(S;\ZZ))\oplus \ZZ\xi$, see~\cite{beau}.
Let $(\cdot,\cdot)$ be the Beauville-Bogomolov bilinear symmetric form on $H^2(S^{[2]})$. It is  
  known~\cite{beau} that
\begin{equation}\label{poldon}
(\mu(\eta),\mu(\eta))=\int_S c_1(\eta)^2,
\quad \mu(H^2(S;\ZZ))\bot \ZZ\xi,
\quad (\xi,\xi)=-2.
\end{equation}
Since $S$ and $S^{[2]}$ are regular varieties we may identify their Picard groups with $H^{1,1}_{\ZZ}(S)$ and $H^{1,1}_{\ZZ}(S^{[2]})$ respectively. Let $C\in \Pic(S)$; abusing notation we will denote by $\mu(C)$ the  class in $\Pic(S^{[2]})$ corresponding to $\mu(\cO_S(C))\in H^{1,1}_{\ZZ}(S)$: if $C$ is an integral curve it is represented by subschemes whose support intersects $C$.
The following is the main result of the present subsection.
\begin{thm}\label{thm:ixhilb}
Let $A\in(\Delta\setminus\Sigma)$ and $[v_0]\in Y_A[3]$ ($=Y_A(3)$ by {Claim 3.11} of~\cite{og5}) - thus $S_A(v_0)$ is a $K3$ surface by~\Ref{crl}{buonogrande}.   Then the following hold:
\begin{itemize}
\item[(1)]
If  $S_A(v_0)$ does not contain lines (true for generic $A$ by~\Ref{prp}{generikappa}) then there exist a choice $\epsilon$ of $\PP^2$-fibration for $X_A$ and an isomorphism. 
\begin{equation}\label{calispera}
\psi\colon {S}_A(v_0)^{[2]}\dashrightarrow X_A^{\epsilon}
\end{equation}
such that 
 \begin{equation}\label{tirodietro}
 \psi^{*}H_A^{\epsilon}\sim 
 \mu(D_A(v_0))-\Delta^{[2]}_{{S}_A(v_0)}.
\end{equation}
\item[(2)]
Let $A$ and $\epsilon$ be arbitrary.  
There exists a bimeromorphic map 
\begin{equation}\label{calimero}
\psi\colon {S}_A(v_0)^{[2]}\dashrightarrow X_A^{\epsilon}
\end{equation}
 such that~\eqref{tirodietro} holds.
\end{itemize}
\end{thm}
\begin{rmk}
Suppose that $S_A(v_0)$ contains a line $L$.  The  restriction of the right-hand side of~\eqref{tirodietro} to $L^{(2)}$ (embedded in $S_A(v_0)^{[2]}$) is $\cO_{L^{(2)}}(-1)$. Since $H_A^{\epsilon}$ is nef we get that in this case Map~\eqref{calimero} cannot be regular.
\end{rmk}
The proof of~\Ref{thm}{ixhilb} will be given after a series of auxiliary results. 
Let $S\subset\PP^6$ be a linearly normal $K3$ surface of genus $6$ such that $\cI_{S/\PP^6}(2)$ is globally generated; then  $S$ is projectively normal and hence Riemann-Roch gives that $\dim |\cI_S(2)|=5$. One defines a rational map ${S}^{[2]}\dashrightarrow
 |\cI_S(2)|^{\vee}$ as follows.  Given $[Z]\in{S}^{[2]}$ we let $\la Z\ra\subset   \PP^5$ be the line spanned by $Z$. We let
\begin{equation}\label{bellaciao}
\begin{matrix}
({S}^{[2]}\setminus \bigcup\limits_
{\text{$L\subset S$
line}}
L^{(2)})
 & \overset{g}{\lra} & |\cI_S(2)|^{\vee}\cong\PP^5 \\
& \\
[Z] & \mapsto & 
\{Q\in |\cI_S(2)|\mid \text{s.t.}\ Q\supset\la Z\ra\}\,.
\end{matrix}
\end{equation}
Let $D$ be a hyperplane divisor on $S$;
one shows (see Claim~(5.16) of~\cite{og1}) that 
\begin{equation}\label{cacciatore}
g^{*}\cO_{\PP^5}(1)\cong \mu({D})-\Delta^{[2]}_{{S}}.
\end{equation}
(Notice that  the set of lines on $S$ is finite and hence $\bigcup_
{L\subset S\text{line}} L^{(2)}$ 
 has codimension $2$ in ${S}^{[2]}$.)
In fact  $g$ can be identified with the map associated to the complete linear system $|(\mu({D})-\Delta^{[2]}_{{S}})|$. We will analyze $g$ under the assumption that $S$ is generic (in a precise sense).
\begin{ass}\label{ass:ipotesse}
Item~(1) of~\Ref{ass}{kappass} holds.  
\begin{equation}
S:=W_K\cap Q
\end{equation}
where $Q\subset\PP(\Ann   K)$ is a quadric intersecting transversely $W_K$.
\end{ass}
Let $S\subset\PP(\Ann   K)$ be as in~\Ref{ass}{ipotesse}. Then $S$ is a linearly normal $K3$ surface of genus $6$ and $\cI_S(2)$ is globally generated. Thus  the map $g$ of~\eqref{bellaciao} is defined. 
Let $F(W_K)$ be the variety parametrizing lines in $W_K$.  Since the set of lines in $S$ is finite (empty for generic $S$ by~\Ref{prp}{generikappa}) 
  we have a  map 
\begin{equation}\label{corda}
\begin{matrix}
(F(W_K)\setminus\{L\mid L\subset S\})& \lra & S^{[2]} \\
L & \mapsto & L\cap Q\,.
\end{matrix}
\end{equation}
\begin{dfn}\label{dfn:pianohilb}
Let $P^0_S\subset S^{[2]}$ be the image of Map~(\ref{corda}) and $P_S$ be its closure in $S^{[2]}$. 
\end{dfn}
We recall that $F(W_K)\cong\PP^2$ by Iskovskih's ~\Ref{prp}{fanoindue}.
\begin{clm}
Let $S\subset\PP(\Ann   K)$ be as in~\Ref{ass}{ipotesse}. Suppose moreover that $S$ contains no lines. Let $C_1,C_2,\ldots,C_s$ be the (smooth) conics contained in $S$ (of course the generic $S$ contains no conics).
Then  $P_S,C_1^{(2)},\ldots,C_s^{[2]}$ are pairwise disjoint subset of $S^{[2]}$. Moreover there exists a biregular morphism
\begin{equation}\label{sgonfio}
c\colon S^{[2]}\lra N(S).
\end{equation}
contracting each of $P_S,C_1^{(2)},\ldots,C_s^{[2]}$.
Thus $N(S)$ is    a compact complex normal space with  
 \begin{equation}\label{singenne}
\sing N(S)=
\{c(P_S),\ldots,c(C^{(2)}),\ldots\mid \text{$C\subset S$  a conic}\}
\end{equation}
and $c$ is an isomorphism of the complement of $P_S\cup C_1^{(2)}\cup\ldots\cup C_s^{[2]}$ onto the smooth locus of $N(S)$. 
The map $g$ (regular on all of $S^{[2]}$ because $S$ contains no lines) descends to a regular map
\begin{equation}
\ov{g}\colon N(S)\to |\cI_S(2)|^{\vee},\qquad \ov{g}\circ c=g\,.
\end{equation}
\end{clm}
\begin{proof}
$P_S$ is isomorphic to $\PP^2$ by Iskovskih's ~\Ref{prp}{fanoindue} and each $C_i^{(2)}$ is  isomorphic to $\PP^2$ because  $C_i$ is a conic. Thus each of $P_S$, $C_i$ can be contracted individually. Let's show that  $P_S,C_1^{(2)},\ldots,C_s^{[2]}$ are pairwise disjoint. Suppose that $[Z]\in P_S\cap C_i^{(2)}$. Let $\Lambda$ be the plane containing $C_i$. Then $\Lambda\cap W_K$ contains the line $\la Z\ra$ and the smooth conic $C_i$. Since $W_K$ is cut out by quadrics it follows that $\Lambda\subset W_K$, that is absurd because $W_K$ contains no planes. This proves that $P_S\cap C_i^{(2)}=\es$. On the other hand there does not exist $[Z]\in C_i^{(2)}\cap C_j^{(2)}$ by~\Ref{crl}{consufan}.  that  $P_S,C_1^{(2)},\ldots,C_s^{[2]}$ are pairwise disjoint. Thus the contraction~\eqref{sgonfio} exists. It remains to prove that $g$ is constant on each of $P_S,C_1^{(2)},\ldots,C_s^{[2]}$. 
In fact if $[Z]\in P_S$ then $g([Z])=|\cI_{W_K}(2)|$, if $[Z]\in C_i^{(2)}$ then 
\begin{equation*}
g([Z])=\{Q\in |\cI_S(2)| \mid Q\supset \la C_i\ra \}.
\end{equation*}
\end{proof}
Now we go back to the \lq\lq general\rq\rq case: we suppose that~\Ref{ass}{ipotesse} holds however $S$ may very well contain lines.
Let
\begin{equation}\label{hilbstella}
S^{[2]}_{\star}:=S^{[2]}\setminus P_S 
\setminus\bigcup\limits_{\text{$R\subset S$ line or conic}}\Hilb^2 R\,.
\end{equation}
(Notice that if $R\subset S$ is a conic which is not smooth then we delete all $[Z]\in S^{[2]}$ such that $Z$ is contained in the scheme $R$.) The following result is essentially {\bf Lemma~3.7} of~\cite{og4}.
\begin{prp}\label{prp:biunivoca}
Suppose that~\Ref{ass}{ipotesse} holds. 
\begin{itemize}
\item[(1)]
The fibers of $g|_{S^{[2]}_{\star}}$ are finite  of cardinality at most $2$ and  the generic fiber has cardinality $2$. 
\item[(2)]
 There exist an open dense subset $\cA\subset S^{[2]}_{\star}$
and  an anti-symplectic (and hence non-trivial)    involution  $\phi\colon \cA\to \cA$ such that 
\begin{equation}\label{invariante}
(g|_{\cA})\circ\phi=g|_{\cA}\,. 
\end{equation}
The induced map
\begin{equation}\label{comequoz}
 \cA/\la \phi\ra\lra 
g(\cA)
\end{equation}
 is a bijection.
\item[(3)]
If in addition   $S$ does not contain lines  $\phi$ descends to a regular involution $\ov{\phi}\colon N(S)\to N(S)$ such that $\ov{g}\circ\ov{\phi}=\ov{g}$ and 
 the induced map
\begin{equation}\label{quozdisc}
j\colon N(S)/\la \ov{\phi}\ra \lra 
g(S^{[2]})
\end{equation}
 is a bijection. Moreover
 \begin{equation}\label{cofisso}
\cod(\Fix(\ov{\phi}),N(S))\ge 2
\end{equation}
where $\Fix(\ov{\phi})$ is the fixed-locus of $\ov{\phi}$.    
\end{itemize}
\end{prp}
 Let $A$ and $[v_0]$ be as in the statement of~\Ref{thm}{ixhilb}: we will perform the 
 key computation one needs to prove that theorem.
 Let $V_0\subset V$ be a codimension-$1$ subspace transversal to $[v_0]$ and such that $\bigwedge^3 V_0\cap A=\{0\}$. Let $\cD$ be Decomposition $V=[v_0]\oplus V_0$ and $S_A^{\cD}$ be given by~(\ref{essdec}) - thus $S_A^{\cD}$ sits in $\PP(\Ann   K^{\cD}_A)\cap\Gr(3,V_0)$ and is isomorphic to $S_A(v_0)$.  
Let $f\in V_0^{\vee}$; we let $q_f$ be the quadratic form on $\bigwedge^3 V_0$ defined by setting
\begin{equation}
q_f(\omega):=\vol_0((f \contr\omega)\wedge\omega)
\end{equation}
where $\vol_0$ is a volume-form on $V_0$.
Then $q_f$ is a Pl\"ucker quadric, in fact we have an isomorphism
\begin{equation}\label{puntura}
\begin{matrix}
V_0^{\vee} & \overset{\sim}{\lra} & H^0(\cI_{\Gr(3,V_0)}(2)) \\
f & \mapsto & q_f\,.
\end{matrix}
\end{equation}
Let $V^{\vee}=[v_0^{\vee}]\oplus V_0^{\vee}$  be the dual decomposition of $\cD$; thus    $v_0^{\vee}\in \Ann   V_0$ and  $v_0^{\vee}(v_0)=1$. 
We have an isomorphism
 \begin{equation}\label{padre}
\begin{matrix}
 [v_0^{\vee}]\oplus V_0^{\vee} & \overset{\sim}{\lra} & H^0(\cI_{S_A^{\vee}}(2))   \\
x v_0^{\vee}+f & \mapsto &  x (r_A^{\cD})^{\vee}+q_f\,.
\end{matrix}
\end{equation}
We let
 \begin{equation}\label{proiso}
\iota\colon |\cI_{S_A^{\cD}}(2)|^{\vee}\overset{\sim}{\lra}\PP(V)
\end{equation}
be the projectivization of the  transpose of~(\ref{padre}). 
\begin{prp}\label{prp:belcalcolo}
Let $A$ and $[v_0]$ be as in the statement of~\Ref{thm}{ixhilb} and keep notation  as above.  Let $g$ be Map~\eqref{bellaciao} for $S_A^{\cD}$ - this makes sense  by~\Ref{crl}{buonogrande}. 
 Then $\iota(\im g)\subset Y_A$. 
\end{prp}
\begin{proof}
Let
\begin{equation}\label{apertodenso}
[Z]\in ((S_A^{\cD})_{\star}^{[2]}\setminus\Delta_{S_A^{\cD}}^{[2]}\setminus 
P_{S_A^{\cD}})\,.
\end{equation}
We will prove that
\begin{equation}
\iota(g([Z])\in Y_A\,.
\end{equation}
This will suffice to prove the lemma because the right-hand side of~(\ref{apertodenso}) is dense in $(S_A^{\cD})_{\star}^{[2]}$ and $Y_A$ is closed. By hypothesis  $Z$ is reduced; thus $Z=\{[\beta],[\beta']\}$ where $\beta,\beta'\in\bigwedge^3 V_0$ are decomposable.  The line $\la[\beta],\beta']\ra$ spanned by $[\beta]$ and $[\beta']$ is not contained in $F^{\cD}_A$ because $[Z]\notin P_{S_A^{\cD}}$. Thus $\la[\beta],\beta']\ra$ is not contained in
$\Gr(3,V_0)$ and it  follows that the vector sub-spaces of $V_0$  supporting the decomposable vectors $\beta$ and $\beta'$ intersect in a $1$-dimensional subspace. Thus there exists a basis $\{v_1,\ldots,v_5\}$ of $V_0$ such that
\begin{equation}
\beta=v_1\wedge v_2\wedge v_3,\quad 
\beta'=v_1\wedge v_4\wedge v_5\,.
\end{equation}
We may assume moreover that $\vol_0(v_1\wedge v_2\wedge v_3\wedge v_4\wedge v_5)=1$. By~\eqref{killk} and~\eqref{errepol}  there exist $\alpha,\alpha'\in\bigwedge^2 V_0$ such that
\begin{equation}
v_0\wedge\alpha+\beta,\ v_0\wedge\alpha'+\beta'\in A,\qquad
\alpha\wedge\beta=\alpha'\wedge\beta'=0\,.
\end{equation}
Since $A$ is Lagrangian we get that
\begin{equation}
\vol_0(\alpha\wedge\beta')=\vol_0(\alpha'\wedge\beta)=:c\,.
\end{equation}
Let $t_0,\ldots,t_5\in\CC$; a straightforward computation gives that
\begin{equation}
(t_0 (r_A^{\cD})^{\vee}+\sum_{i=1}^5 t_i q_{v_i^{\vee}})(\beta+\beta')=
2ct_0+2t_1\,.
\end{equation}
Thus 
\begin{equation}\label{imagogi}
\iota(g([Z]))=[cv_0+v_1]\,.
\end{equation}
It remains to prove that
\begin{equation}\label{appartiene}
[cv_0+v_1]\in Y_A\,.
\end{equation}
Let $K^{\cD}_A$ be as in~(\ref{kappagrasso}); we claim that it suffices to prove that there exist $(x,x')\in(\CC^2\setminus\{(0,0)\})$  and $\kappa\in K^{\cD}_A$ such that
\begin{equation}\label{desiderata}
(cv_0+v_1)\wedge
(x(v_0\wedge\alpha+\beta)+x'(v_0\wedge\alpha'+\beta')+
v_0\wedge\kappa)=0\,.
\end{equation}
In fact assume that~(\ref{desiderata}) holds. Then 
\begin{equation}
0\not=(x(v_0\wedge\alpha+\beta)+x'(v_0\wedge\alpha'+\beta')+
v_0\wedge\kappa)\in A\cap F_{cv_0+v_1}.
\end{equation}
(The inequality holds  because $\beta,\beta'$ are linearly independent.) A straightforward computation gives that~(\ref{desiderata}) is equivalent to
\begin{equation}\label{olimpia}
x(c\beta-v_1\wedge\alpha)+x'(c\beta'-v_1\wedge\alpha')=
v_1\wedge\kappa\,.
\end{equation}
As is easily checked we have
\begin{equation}\label{podio}
(c\beta-v_1\wedge\alpha),\ (c\beta'-v_1\wedge\alpha')\in
([v_1]\wedge(\bigwedge^ 2\la v_2,v_3,v_4,v_5\ra))\cap
\{v_2\wedge v_3,\ v_4\wedge v_5\}^{\bot}
\end{equation}
where perpendicularity is with respect to wedge-product followed by $\vol_0$. 
Multiplication by $v_1$ gives an injection $K^{\cD}_A\hra
([v_1]\wedge(\bigwedge^ 2\la v_2,v_3,v_4,v_5\ra))$; in fact   no non-zero element of $K^{\cD}_A$ is decomposable  because $A\notin\Sigma$. Since the right-hand side of~(\ref{podio}) has dimension $4$ and $\dim K^{\cD}_A=3$ we get that there exists $(x,x')\in(\CC^2\setminus\{(0,0)\})$ such that~(\ref{olimpia}) holds. 
\end{proof}
\begin{lmm}\label{lmm:rivcon}
Let $A\in(\lagr\setminus\Sigma)$.  Then $Y_A(1)$ is not empty, the topological double cover $f_A^{-1}Y_A(1)\to Y_A(1)$ is not trivial and $Y_A$ is integral.
\end{lmm}
\begin{proof}
By~\Ref{clm}{nodeco} we know that $Y_A[3]$ is finite. On the other hand $(Y_A[2]\setminus Y_A[3])$ is a smooth surface   - see Proposition~2.8 of~\cite{og2}. Since $\sing Y_A\subset Y_A[2]$ it follows that $Y_A$ is integral and $Y_A(1)$ is connected. Let $[v_0]\in (Y_A[2]\setminus Y_A[3])$. 
By~\Ref{prp}{gradodue} we know that $f_A^{-1}([v_0])$ is a singleton $\{q\}$. Moreover $X_A$ is smooth at $q$ by~\Ref{lmm}{critliscio}. Thus there exists an open neighborhood $U$ of $[v_0]$ in $Y_A$ such that $f_A^{-1}U$ is smooth. Moreover $(f_A^{-1}Y_A[2])\cap f_A^{-1} U$ is nowhere dense in $f_A^{-1} U$. Since $f_A^{-1} U$ is smooth the complement $f_A^{-1}(Y_A(1)\cap U)$ is connected.  Since $Y_A(1)$ is connected it follows that 
 $f_A^{-1}Y_A(1)$ is connected.
\end{proof}
\begin{prp}
Keep hypotheses and notation as in~\Ref{prp}{belcalcolo}. Then $\iota(\ov{\im g})=Y_A$. 
\end{prp}
\begin{proof}
By Item~(1) of~\Ref{prp}{biunivoca} the map $g$ has finite generic fiber and hence $\dim\ov{\im g}=4$. By~\Ref{prp}{belcalcolo} we get that $\iota(\ov{\im   g})$ is an irreducible component of $Y_A$. On the other hand   $Y_A$ is irreducible by~\Ref{lmm}{rivcon};  it follows that $\iota(\ov{\im g})=Y_A$.
\end{proof}
\begin{rmk}\label{rmk:puntovuz}
Keep notation as in~\Ref{prp}{belcalcolo};  then
\begin{equation}
\iota\circ g(P^0_{S_A^{\cD}})=\iota(H^0(\cI_{F^{\cD}_A}(2)))=[v_0].
\end{equation}
\end{rmk}
\noindent
{\it Proof of~\Ref{thm}{ixhilb}.} 
Let's prove that Item~(1) holds.  Let $A$ and $[v_0]$ be as in the statement of~\Ref{thm}{ixhilb}. Let $V_0\subset V$ be a codimension-$1$ subspace transversal to $[v_0]$ and such that $\bigwedge^3 V_0\cap A=\{0\}$. Let $\cD$ be Decomposition $V=[v_0]\oplus V_0$. In order to simplify notation we set $S=S_A^{\cD}$; thus $S\cong S_A(v_0)$ and by hypothesis $S$ does not contain lines. 
Let $j$ be the map of~\eqref{quozdisc}; by~\Ref{prp}{belcalcolo} the composition $\iota\circ j$ is a map
\begin{equation}\label{pettirosso}
\iota\circ j\colon N(S)/\la\ov{\phi}\ra \lra Y_A\,.
\end{equation}
We claim that $\iota\circ j$ is an isomorphism: 
in fact it has finite fibers and is birational by~\Ref{prp}{biunivoca}, since $\dim \sing Y_A=2$ (because $A\notin\Sigma$)  the hypersurface $Y_A$ is normal and hence  $\iota\circ j$ is an isomorphism. Let $\pi\colon N(S)\to N(S)/\la\ov{\phi}\ra$ be the quotient map. 
By~\eqref{cofisso} the singular locus of $N(S)/\la\ov{\phi}\ra$ is the image of $\Fix(\ov{\phi})$ (and thus isomorphic to $\Fix(\ov{\phi})$); since~(\ref{pettirosso}) is an isomorphism  we get that
\begin{equation}\label{rivestodoppio}
\begin{matrix}
N(S)\setminus \Fix(\ov{\phi}) & \lra & Y_A^{sm} \\
x & \mapsto & \iota\circ j\circ\pi(x)
\end{matrix}
\end{equation}
 is a topological covering of degree $2$. 
We claim that
\begin{equation}\label{zetadue}
\pi_1(Y_A^{sm})\cong\ZZ/(2)\,.
\end{equation}
In fact $(N(S)\setminus \Fix(\ov{\phi}))\cong(S^{[2]}\setminus(P_S\cup \Fix(\phi|_{S^{[2]}\setminus P_S}))$. Since  $(P_S\cup \Fix(\phi|_{S^{[2]}\setminus P_S}))$ is of codimension $2$ in the simply connected
 manifold
 $S^{[2]}$ we get that $(N(S)\setminus \Fix(\ov{\phi}))$ is simply connected. Thus~(\ref{rivestodoppio}) is the universal covering of $Y_A^{sm}$ and we get~(\ref{zetadue}). On the other hand $Y_A^{sm}\subset Y_A(1)$ by {Corollary~1.5} of~\cite{og5} and thus by~\Ref{lmm}{rivcon} we get that $f_A^{-1}Y_A^{sm}\to Y_A^{sm}$ is the universal covering of $Y_A^{sm}$ as well.
 Hence both $X_A$ and $N(S)$ are normal completions of the universal cover of $Y_A^{sm}$ such that the extended maps to $Y_A$ are finite; it follows that they are isomorphic (over $Y_A$).
  The singular locus of $N(S)$ is given by~\eqref{singenne}. On the other hand $\sing X_A=Y_A[3]$. By~\Ref{rmk}{puntovuz} we can order the set of (smooth) conics on $S$, say $C_1,\ldots,C_s$ and the set of points in $Y_A[3]$ different from $[v_0]$, say $[v_1],\ldots,[v_s]$ so that 
\begin{equation}\label{buonatale}
\ov{\psi}(c(P_S))=[v_0],\qquad\ov{\psi}(c(C_i^{(2)}))=[v_i],\quad
1\le i\le s.
\end{equation}
(Recall~\Ref{rmk}{puntovuz}.) Let $\epsilon_0$ be a choice of $\PP^2$-fibration for $X_A$; then $\ov{\psi}$ defines a birational map $\psi_0\colon S^{[2]}\dra X_A^{\epsilon_0}$ such that
\begin{equation}\label{wurstel}
\psi_0^{*}H_A^{\epsilon_0}\cong\mu(D)-\Delta_S^{[2]}
\end{equation}
where $D$ is the hyperplane class of $S$ (thus $(S,D)$ is isomorphic to $(S_A(v_0),D_A(v_0))$). The birational map $\psi_0$ is an isomorphism away from 
\begin{equation}\label{undisg}
P_S\cup C_1^{(2)}\cup\ldots\cup C_s^{(2)}.
\end{equation}
It follows that $\psi_0$ is the flop of a collection of irreducible components of~(\ref{undisg}). By~\Ref{prp}{fuorisig} we get that there exists   a choice of $\PP^2$-fibration for $X_A$, call it $\epsilon$, such that the corresponding birational map $\psi\colon S^{[2]}\dra X_A^{\epsilon}$ is biregular. Equation~(\ref{tirodietro}) follows from~(\ref{wurstel}). This finishes the proof  that Item~(1) holds. 
Item~(2) follows from Item~(1) and a specialization argument - we leave the details to the reader.  
\qed

\vskip 2mm
\n
We close the present subsection by reproving a result of ours. Let $h_A:=c_1(\cO_{X_A}(H_A))$. 
\begin{thm}[{\rm O'Grady~\cite{og2}}]\label{thm:beatrice}
Let $A\in\lagr^0$. Then $X_A$ is a deformation of $(K3)^{[2]}$ and $(h_A,h_A)_{X_A}=2$. Any small deformation of $(X_A,H_A)$ (i.e.~a small deformation of $X_A$ keeping $h_A$ of type $(1,1)$) is isomorphic to $(X_B,H_B)$ for some $B\in\lagr^0$. 
\end{thm}
\begin{proof}
Let $A_0\in (\Delta\setminus\Sigma)$ and $[v_0]\in Y_{A_0}[3]$.  Suppose moreover that  $S_{A_0}(v_0)$ does not contain lines. By~\Ref{thm}{ixhilb} there exists a choice $\epsilon$ of $\PP^2$-fibration for $X_{A_0}$ such that  
we have an isomorphism 
\begin{equation}
\psi\colon S^{[2]}\overset{\sim}{\lra} X_{A_0}^{\epsilon},
\qquad
\psi^{*}H_{A_0}^{\epsilon}\sim\mu(D_A(v_0))-\Delta_{S_{A_0}(v_0)}^{[2]}.
\end{equation}
On the other hand  $(X_A,H_A)$ is a deformation of $(X_{A_0}^{\epsilon},H_{A_0}^{\epsilon})$ by~\Ref{crl}{unadef}.
This proves that $(X_A,H_A)$ is a deformation of $(S^{[2]},(\mu(D_A(v_0))-\Delta_{S_{A_0}(v_0)}^{[2]}))$.
By~\eqref{poldon} we get that $(h_A,h_A)_{X_A}=2$. Lastly we prove that an arbitrary small deformation of $(X_A,H_A)$ is isomorphic to  $(X_{A'},H_{A'})$ for some $A'\in\lagr^0$. The deformation space of $(X_A,H_A)$ has dimension given by
\begin{equation}\label{defixacca}
\dim \Def(X_A,H_A)=h^{1,1}(X_A)-1=20\,.
\end{equation}
On the other hand $\lagr^0$ is contained in the locus of points in $\lag$ which are stable for the natural (linearized) $PL(V)$-action - this is proved in~\cite{og2}. Thus  by varying $A\in\lagr$ we get
\begin{equation}\label{defepw}
\dim \lagr-\dim SL(V)=55-35=20
\end{equation}
moduli of double EPW-sextics. Since~(\ref{defixacca}) and~(\ref{defepw}) are equal we conclude that an arbitrary small deformation of $(X_A,H_A)$ is isomorphic to  $(X_{B},H_{B})$ for some $B\in\lagr^0$.
\end{proof}
\section{Appendix: Three-dimensional sections of $\Gr(3,\CC^5)$}\label{sec:wedding}
\setcounter{equation}{0}
In the present section $V_0$ is a complex vector-space of dimension $5$. 
Choose a volume form $\vol_0$ on $V_0$; it defines an isomorphism
\begin{equation}\label{zeppa}
\begin{matrix}
\bigwedge^ 2 V_0 & \overset{\sim}{\lra} & \bigwedge^ 3 V^{\vee}_0 \\
\alpha & \mapsto & \omega\mapsto \vol_0(\alpha\wedge\omega)
\end{matrix}
\end{equation}
Let $K\subset\bigwedge^ 2 V_0$ be a $3$-dimensional subspace 
  such that either 
\begin{equation}\label{nointer}
\PP(K)\cap \Gr(2,V_0)=\es
\end{equation}
or else
\begin{equation}\label{puntosolo}
\PP(K)\cap \Gr(2,V_0)=\{[\kappa_0]\}= 
\PP(K)\cap T_{[\kappa_0]}\Gr(2,V_0)\,.
\end{equation}
In other words either $\PP(K)$ does not intersects $\Gr(2,V_0)$ or else the scheme-theoretic intersection is a single reduced point.
We will describe
\begin{equation}
W_{K}:=\PP(\Ann   K)\cap \Gr(3,V_0)
\end{equation}
 First we recall that the dual of $\Gr(3,V_0)$ is $\Gr(2,V_0)$. More precisely let $[\alpha]\in\PP(\bigwedge^ 2 V_0)$: then
\begin{equation}\label{sezsing}
\sing(\PP(\Ann   \alpha)\cap\Gr(3,V_0))=
\{U\in\Gr(3,V_0)\mid U\supset \supp \alpha\}.
\end{equation}
In particular  $\PP(\Ann   \alpha)$ is tangent to $\Gr(3,V_0)$ if and only if $[\alpha]\in\Gr(2,V_0)$ (and if that is the case it is tangent along a $\PP^2$).
Secondly we record the following observation (the proof is an easy exercise).
\begin{lmm}\label{lmm:senullo}
Let $U\subset V_0$ be a codimension-$1$ subspace. Let $\alpha\in\bigwedge^ 2 V_0$. Then
\begin{equation}
\alpha\wedge(\bigwedge^ 3 U)=0
\end{equation}
if and only if $ \supp \alpha\subset U$. 
\end{lmm} 
We recall the following result of Iskovskih.
\begin{prp}[{\rm Iskovskih~\cite{iskovskih}}]\label{prp:fanoindue}
Keep notation as above. Let $K\subset\bigwedge^ 2 V_0$ be a $3$-dimensional subspace 
  such that~(\ref{nointer}) holds. Then 
 \begin{itemize}
\item[(1)]
$W_{K}$ is a smooth Fano $3$-fold of degree $5$ with $\omega_{W_{K}}\cong\cO_{W_{K}}(-2)$,
\item[(2)]
the Fano variety $F(W_{K})$ parametrizing lines on $W_{K}$ (reduced structure) is isomorphic to $\PP^2$,
\item[(3)]
the projective equivalence class of $W_{K}$ does not depend on $K$.
\end{itemize}
\end{prp}
\begin{prp}\label{prp:singfano}
Keep notation as above. Let $K\subset\bigwedge^ 2 V_0$ be a sub vector-space 
 of dimension $3$ such that~\eqref{puntosolo} holds. Then 
 $W_{K}$ is a singular Fano  $3$-fold of degree $5$ with $\omega_{W_{K}}\cong\cO_{W_{K}}(-2)$ and one singular point which is ordinary quadratic and belongs to
\begin{equation}
\{U\in\Gr(3,V_0)\mid U\supset \supp\kappa_0\}.
\end{equation}
\end{prp}
 \begin{proof}
If $\kappa\in(K\setminus[\kappa_0])$ then $\kappa$ is not decomposable and hence $\PP(\Ann   \kappa)$ is transverse to $\Gr(3,V_0)$; by~(\ref{sezsing}) we get that 
\begin{equation}\label{nelpiano}
\sing W_{K}=
\{U\in\Gr(3,V_0)\mid U\supset \supp\kappa_0\}
\cap\PP(\Ann   K)\,.
\end{equation}
We claim that the above intersection consists of one point. First notice that we have a natural identification 
\begin{equation}\label{miomartello}
\{U\in\Gr(3,V_0)\mid U\supset \supp\kappa_0\}\cong
\PP(V_0/ \supp\kappa_0)
\end{equation}
and a linear map
\begin{equation}
\begin{matrix}
K & \overset{\nu}{\lra} & (V_0/ \supp\kappa_0)^{\vee} \\
\kappa & \mapsto & (\ov{v}\mapsto \vol_0(v\wedge\kappa_0\wedge\kappa))
\end{matrix}
\end{equation}
where $v\in V_0$ and $\ov{v}$ is its class in $V_0/ \supp\kappa_0$.
Given~(\ref{nelpiano}) and Identification~(\ref{miomartello}) we get that
\begin{equation}
\sing W_{K}=\PP(\Ann   \im \nu)\,.
\end{equation}
 Of course $\kappa_0\in\ker\nu$ and hence
in order to prove that $\sing W_{K}$ is a singleton it suffices to prove that $\ker\nu=[\kappa_0]$. If $\kappa\in(K\setminus[\kappa_0])$ then $\kappa_0\wedge\kappa\not=0$; in fact this follows from~(\ref{puntosolo}) together with the equality 
\begin{equation}
\PP\{\kappa\in\bigwedge^ 2 V_0\mid 
 \kappa_0\wedge\kappa=0\}=T_{[\kappa_0]}\Gr(2,V_0)\,.
\end{equation}
Since $\kappa_0\wedge\kappa\not=0$ we have  $\nu(\kappa)\not=0$. 
This proves that $\sing W_{K}$ consists of a single point. The formula for the dualizing sheaf of $W_{K}$ follows at once from adjunction. It remains to prove that $W_{K}$ has a single singular point and that it is an ordinary quadratic  point. 
Let $\wt{W}_{K}\subset\PP(\supp\kappa_0)\times
\PP(V_0/ \supp\kappa_0)\times W_{K}$ be the closed subset defined by
\begin{equation}
\wt{W}_{K}:=\{([v],U,W)\mid
v\in W\subset U\}\,.
\end{equation}
The projection $ \wt{W}_{K}\to \PP(V_0/ \supp\kappa_0)$ is a $\PP^1$-fibration and hence $ \wt{W}_{K}$ is smooth. One shows that the projection $\pi\colon\wt{W}_{K}\to W_{K}$ is the blow-up of $\sing W_{K}$. Moreover $\pi^{-1}(\sing W_{K})\cong\PP^1\times\PP^1$ and one gets that the singularity of $W_{K}$ is ordinary quadratic. 
\end{proof}
Our last result is about the base-locus of $3$-dimensional linear systems of quadrics containing $W_K$ for $K\subset\bigwedge^ 2 V_0$ a $3$-dimensional subspace such that~(\ref{nointer}) holds. First we consider the analogous question for the Grassmannian $\Gr(3,\bigwedge^ 3 V_0)$. Let's consider the rational map 
\begin{equation}
\PP(\bigwedge^ 3 V_0)\overset{\Phi}{\dashrightarrow}
|\cI_{\Gr(3,V_0)}(2)|^{\vee}\cong\PP(V_0)
\end{equation}
where the last isomorphism is given by~(\ref{puntura}). 
Let $Z\subset \PP(\bigwedge^ 3 V_0)\times\PP(V_0)$ be the incidence subvariety defined by
\begin{equation}
Z:=\{([\omega],[v])\mid v\wedge\omega=0\}\,.
\end{equation}
Then we have a commutative triangle 
\begin{equation}\label{ranieri}
\xymatrix{  
& Z \ar^{\wt{\Phi}}[dr] \ar_{\Psi}[dl]  & \\
 \PP(\bigwedge^ 3 V_0) &\overset{\Phi}{\dashrightarrow} & \PP( V_0)}
\end{equation}
where $\Psi$ and $\wt{\Phi}$ are the restrictions to $Z$ of the two projections of $\PP(\bigwedge^ 3 V_0)\times\PP(V_0)$. Moreover $\Psi$ is the blow-up of $\Gr(3,V_0)$. In particular the following holds: if $\omega\in\bigwedge^ 3 V_0$ is not decomposable then there  exists a unique $[v]\in\PP(V_0)$ such that $v\wedge\omega=0$ and moreover $\Phi([\omega])=[v]$. Let $[v]\in\PP(V_0)$; by~(\ref{puntura}) we may view $\Ann  (v)\subset V_0^{\vee}$ as a hyperplane in $|\cI_{\Gr(3,V_0)}(2)|$; by commutativity of~(\ref{ranieri}) we have
\begin{equation}\label{intspezzata}
\bigcap\limits_{f\in \Ann  (v)}V(q_f)=\Gr(3,V_0)
\cup\{[\omega]\in\PP(\bigwedge^ 3 V_0)\mid v\wedge\omega=0\}.
\end{equation}
\begin{prp}\label{prp:baseweb}
Let $K\subset\bigwedge^ 2 V_0$ be a $3$-dimensional subspace such that~(\ref{nointer}) holds. Let $L\subset|\cI_{W_K}(2)|$ be a hyperplane (here $\cI_{W_K}$ is the ideal sheaf of $W_K$ in $\PP(\Ann   K)$). Then
\begin{equation}
\bigcap\limits_{t\in L}Q_t=W_K\cup R_L
\end{equation}
where $R_L$ is a plane. Moreover $W_K\cap R_L$ is a  conic.
\end{prp}
\begin{proof}
Restriction to $\PP(\Ann   K)$ defines an isomorphism
\begin{equation}
|\cI_{\Gr(3,V_0)}(2)|\overset{\sim}{\lra} |\cI_{W_K}(2)|\,.
\end{equation}
By~(\ref{puntura}) we get that we may identify $L$ with $\PP(\Ann  (v))$ for a well-defined $[v]\in\PP(V_0)$ and each quadric $Q_t$ for $t\in L$ with $\PP(\Ann   K)\cap V(q_f)$ for a suitable $[f]\in \PP(\Ann  (v))$. By~(\ref{intspezzata}) we have
\begin{equation}
\bigcap\limits_{f\in \Ann  (v)}(\PP(\Ann   K)\cap V(q_f))=
W_K\cup R_L
\end{equation}
where
\begin{equation}
R_L:=\PP(\Ann   K)\cap 
\{[\omega]\in\PP(\bigwedge^ 3 V_0)\mid v\wedge\omega=0\}.
\end{equation}
 Thus $R_L$ is a linear space of dimension at least $2$. Now notice that we have an isomorphism
 \begin{equation}
\begin{matrix}
\bigwedge^ 2 (V_0/[v]) & \overset{\sim}{\lra} & \{[\omega]\in\PP(\bigwedge^ 3 V_0)\mid v\wedge\omega=0\} \\
\ov{\alpha} & \mapsto & v\wedge\alpha
\end{matrix}
\end{equation}
where $\alpha\in \bigwedge^ 2 V_0$ is an element mapped to $\ov{\alpha}$ by the quotient map $\bigwedge^ 2 V_0\to \bigwedge^ 2 (V_0/[v])$. Since $\dim(V_0/[v])=4$ the Grassmannian $\Gr(2,V_0/[v])$ is a quadric hypersurface in $\PP(\bigwedge^ 2(V_0/[v]))$; it follows that either $R_L\subset W_K$ or $R_L\cap W_K$ is a quadric hypersurface in $R_L$. By Lefschetz $\Pic(W_K)$ is generated by the hyperplane class; it follows that $W_K$ contains no planes and no  quadric surfaces. Thus necessarily $\dim R_L=2$, moreover $R_L\not\subset W_K$ and the intersection $R_L\cap W_K$ is a conic.
\end{proof}
\begin{crl}\label{crl:consufan}
Let $K\subset\bigwedge^ 2 V_0$ be a $3$-dimensional subspace such that~(\ref{nointer}) holds and  $\cC(W_K)$ be  the variety parametrizing conics on $W_K$ (reduced structure). Then  we have an isomorphism
\begin{equation}
\begin{matrix}
|\cI_{W_K}(2)|^{\vee} & \overset{\sim}{\lra} & \cC(W_K) \\
L & \mapsto & R_L\cap W_K
\end{matrix}
\end{equation}
where $R_L$ is as in~\Ref{prp}{baseweb}. Moreover given $Z\in W_K^{[2]}$ there exists a unique conic containing $Z$ namely $R_L\cap W_K$ where $L\in |\cI_{W_K}(2)|^{\vee}$ is the hypeprlane of quadrics containing $\la Z\ra$. 
\end{crl}
\end{document}